\documentclass[a4,12pt,reqno]{amsart}
\usepackage{amsmath}
\usepackage{amsthm}
\usepackage{graphicx}
\usepackage{amsfonts}
\usepackage{bm}
\usepackage{amssymb}
\usepackage{enumerate}
\usepackage[margin=1.1in]{geometry}
\usepackage{xcolor}

\numberwithin{equation}{section}

\newtheorem{theorem}{Theorem}[section]
\newtheorem{lemma}[theorem]{Lemma}
\newtheorem{proposition}[theorem]{Proposition}

\theoremstyle{definition}

\newtheorem{innercustomthm}{Assumption}
\newenvironment{assumption}[1]
  {\renewcommand\theinnercustomthm{#1}\innercustomthm}
  {\endinnercustomthm}

\def\E{{\mathbb E}}
\def\R{{\mathbb R}}

\def\N{{\mathbb N}}

\def\PP{{\mathbb P}}
\def\FF{{\mathbb F}}

\def\F{{\mathcal F}}

\begin{document}
\title[Inverting the Markovian projection]{Inverting the Markovian projection, with an application to local stochastic volatility models}
\author{Daniel Lacker, Mykhaylo Shkolnikov, and Jiacheng Zhang}
\address{IEOR Department, Columbia University, New York, NY 10027}
\email{daniel.lacker@columbia.edu}
\address{ORFE Department, Bendheim Center for Finance, and Program in Applied \& Computational Mathematics, Princeton University, Princeton, NJ 08544.}
\email{mshkolni@gmail.com}
\address{ORFE Department, Princeton University, Princeton, NJ 08544.}
\email{jiacheng@princeton.edu}
\footnotetext[1]{M.~Shkolnikov is partially supported by the NSF grant DMS-1811723 and a Princeton SEAS innovation research grant.}

\begin{abstract}
We study two-dimensional stochastic differential equations (SDEs) of McKean--Vlasov type in which the conditional distribution of the second component of the solution given the first enters the equation for the first component of the solution. Such SDEs arise when one tries to invert the Markovian projection developed in \cite{gyongy1986mimicking}, typically to produce an It\^o process with the fixed-time marginal distributions of a given one-dimensional diffusion but richer dynamical features. We prove the strong existence of stationary solutions for these SDEs, as well as their strong uniqueness in an important special case. Variants of the SDEs discussed in this paper enjoy frequent application in the calibration of local stochastic volatility models in finance, despite the very limited theoretical understanding.   
\end{abstract}

\maketitle


\section{Introduction}

We consider a class of stochastic differential equations (SDEs) that arises naturally when one attempts to invert the \emph{Markovian projection}, a concept originating from a celebrated theorem of \textsc{Gy\"ongy} \cite[Theorem 4.6]{gyongy1986mimicking}. The idea of a Markovian projection, often also attributed to \textsc{Krylov} \cite{krylov1984relation}, lies in finding a diffusion which ``mimicks'' the fixed-time marginal distributions of an It\^o process. We quote here a version due to \textsc{Brunick} and \textsc{Shreve} \cite[Corollary 3.7]{brunick2013mimicking}, which significantly relaxes the assumptions on the coefficients in \cite{gyongy1986mimicking}.

\begin{proposition}[Markovian projection, \cite{brunick2013mimicking}] \label{th:intro:Markovprojection}
Let $(b_t)_{t \ge 0}$ and $(\sigma_t)_{t \ge 0}$ be adapted real-valued processes defined on a stochastic basis $(\Omega,\F,\FF,\PP)$ supporting an $\FF$-Wiener process $(W_t)_{t \ge 0}$ and such that $\E\big[\int_0^t |b_s| + \sigma_s^2\,\mathrm{d}s\big] < \infty$ for each $t > 0$. Suppose a process $(X_t)_{t \ge 0}$ satisfies
\[
\mathrm{d}X_t = b_t\,\mathrm{d}t + \sigma_t\,\mathrm{d}W_t.
\]
Then there are measurable functions $\widehat b:\, [0,\infty) \times \R \to \R$ and $\widehat \sigma :\, [0,\infty) \times \R \to \R$ so that:
\begin{enumerate}
\item[{\rm (i)}] For a.e.\ $t \ge 0$, one has the a.s.\ equalities 
\[
\widehat b(t,X_t) = \E[b_t  |  X_t] \quad\text{and}\quad \widehat \sigma(t,X_t)^2 = \E[\sigma_t^2|X_t].
\]
\item[{\rm (ii)}] There exists a weak solution of the SDE
\begin{align}
\mathrm{d}\widehat{X}_t = \widehat b(t,\widehat{X}_t)\,\mathrm{d}t + \widehat \sigma(t,\widehat{X}_t)\,\mathrm{d}W_t \label{intro:mimicSDE}
\end{align}
with the property that $\widehat{X}_t\stackrel{d}{=}X_t$ for all $t\ge0$, where $\stackrel{d}{=}$ denotes equality in law.
\end{enumerate}
\end{proposition}

\smallskip

Herein, we are interested in inverting the Markovian projection, that is, in finding a different It\^o process with the fixed-time marginal distributions matching those of a given one-dimensional diffusion. This problem appears, for example, in the calibration procedure for local stochastic volatility models in finance (see  \cite{lipton2002masterclass}, \cite{liptonmcghee2002masterclass}, \cite{piterbarg2006markovian}, \cite{guyon2011solved}, \cite{guyon2012being}, \cite[Chapter 11]{guyon2013nonlinear}, \cite{tian2015calibrating}, \cite{abergel2017nonparametric}, \cite{saporito2017calibration}, as well as further below in this introduction). Given a one-dimensional diffusion
\begin{equation} \label{eq:Xhat}
\mathrm{d}\widehat{X}_t=b_1(\widehat{X}_t)\,\mathrm{d}t+\sigma_1(\widehat{X}_t)\,\mathrm{d}\widehat{W}_t,
\end{equation}
Proposition \ref{th:intro:Markovprojection} suggests the ansatz 
\begin{equation}\label{eq:SDE_gen}
\mathrm{d}X_t=\gamma_t\,\frac{b_1(X_t)}{\E[\gamma_t|X_t]}\,\mathrm{d}t + \zeta_t\,\frac{\sigma_1(X_t)}{\sqrt{\E[\zeta_t^2|X_t]}}\,\mathrm{d}W_t
\end{equation}
with adapted processes $(\gamma_t)_{t\ge0}$ and $(\zeta_t)_{t\ge0}$. However, due to the presence of the conditional expectations in the equation, for an arbitrary choice of $(\gamma_t)_{t\ge0}$ and $(\zeta_t)_{t\ge0}$ the construction of an It\^o process $(X_t)_{t\ge0}$ satisfying \eqref{eq:SDE_gen} seems completely out of reach. Therefore, we specialize to the setting where $(\Omega,\F,\FF,\PP)$ supports an $\FF$-Wiener process $(B_t)_{t\ge0}$ independent of $(W_t)_{t\ge0}$. With a one-dimensional diffusion
\begin{equation}\label{eq:Y_def}
\mathrm{d}Y_t=b_2(Y_t)\,\mathrm{d}t + \sigma_2(Y_t)\,\mathrm{d}B_t,
\end{equation}
we further set $\gamma_t=h(Y_t)$ and $\zeta_t=f(Y_t)$, for all $t\ge0$ and some measurable functions $h$ and $f$. As a result, we are led to consider the two-dimensional SDE  
\begin{align}
\begin{cases}
\;\mathrm{d}X_t = b_1(X_t)\,\frac{h(Y_t)}{\E[h(Y_t)|X_t]}\,\mathrm{d}t + \sigma_1(X_t)\,\frac{f(Y_t)}{\sqrt{\E[f^2(Y_t)|X_t]}}\,\mathrm{d}W_t, \\
\;\mathrm{d}Y_t \,=\, b_2(Y_t)\,\mathrm{d}t + \sigma_2(Y_t)\,\mathrm{d}B_t.
\end{cases} \label{intro:mainSDEsystem}
\end{align}

\smallskip

Our first main theorem yields the strong existence of a stationary solution for the SDE \eqref{intro:mainSDEsystem} under the following assumption.

\begin{assumption}{\textbf{A}} \label{asmp:main} The functions $(b_1,b_2,\sigma_1,\sigma_2,h,f)$ are measurable and satisfy:
\begin{enumerate}[(a)]
\item There exist constants $c,C_1,C_2\in(0,\infty)$ such that for $i=1,2$ and all $x\in\R$: 	
\begin{equation}\label{b_asmp}
x b_i(x) \le -cx^2+C_1
\quad\text{and}\quad |b_i(x)|\le C_2(1+|x|). 
\end{equation}
The functions $\sigma_1$ and $\sigma_2$ are bounded above and below by positive constants and possess bounded derivatives, $\sigma'_1$ and $\sigma'_2$. 
\item The functions $h$ and $f$ are bounded above and below by positive constants, and $f$ admits a bounded derivative $f'$.  
\end{enumerate}
\end{assumption}

Before stating our main theorem, we first state a straightforward and standard lemma to fix terminology. Both Lemma \ref{le:intro:onedimensional} and Theorem \ref{th:intro:main} are proven in Sections \ref{se:2.1} and \ref{se:2.2}.

\begin{lemma} \label{le:intro:onedimensional}
Under Assumption \ref{asmp:main}, the one-dimensional SDE \eqref{eq:Xhat} admits a unique in law strong solution starting from any initial position, and there is a unique in law solution satisfying the stationarity property $\widehat{X}_t \stackrel{d}{=} \widehat{X}_0$ for all $t\ge 0$. In addition, the same claims are true for the one-dimensional SDE
\begin{align}
\mathrm{d}\widehat{Y}_t \,=\, b_2(\widehat{Y}_t)\,\mathrm{d}t + \sigma_2(\widehat{Y}_t)\,\mathrm{d}\widehat{B}_t. \label{eq:Yhat}
\end{align}
In stationarity, the laws of each $\widehat{X}_t$ and $\widehat{Y}_t$ admit densities $m_1$ and $m_2$, respectively, where
\begin{align}
m_i(x)\!=\! \frac{1}{\sigma_i^2(x)}\!\exp\bigg(\!\int_0^x\frac{2b_i(a)}{\sigma_i^2(a)}\,\mathrm{d}a\!\bigg) \!\bigg/\! \int_{\R} \frac{1}{\sigma_i^2(a_1)}\,\exp\!\bigg(\!\int_0^{a_1} \frac{2b_i(a_2)}{\sigma_i^2(a_2)}\,\mathrm{d}a_2\!\bigg)\,\mathrm{d}a_1, \; i=1,2. \label{def:marginals}
\end{align}
\end{lemma}

\begin{theorem} \label{th:intro:main}
Under Assumption \ref{asmp:main}, there exists a weak solution $(X_t,Y_t)_{t \ge 0}$ of the SDE \eqref{intro:mainSDEsystem} satisfying the stationarity property $(X_t,Y_t)\stackrel{d}{=}(X_0,Y_0)$ for all $t\ge0$. Moreover, any stationary weak solution $(X_t,Y_t)_{t \ge 0}$ of \eqref{intro:mainSDEsystem} is strong, and the following hold:
\begin{enumerate}
\item[{\rm (i)}] $X_t\stackrel{d}{=}\widehat{X}_t$ for all $t\ge0$, where $(\widehat{X}_t)_{t\ge0}$ is the unique stationary solution of \eqref{eq:Xhat}.
\item[{\rm (ii)}] $Y_t\stackrel{d}{=}\widehat{Y}_t$ for all $t\ge0$, where $(\widehat{Y}_t)_{t\ge0}$ is the unique stationary solution of \eqref{eq:Yhat}.
\item[{\rm (iii)}] The law of $(X_0,Y_0)$ admits a density $p$ with $\int_{\R^2} \frac{|\nabla p|^2}{p}\,\mathrm{d}x\,\mathrm{d}y<\infty$.
\end{enumerate}
\end{theorem}

\smallskip

Parts (i) and (ii) of Theorem \ref{th:intro:main}, which are consequences of Lemma \ref{le:intro:onedimensional} and Proposition \ref{th:intro:Markovprojection}, ensure that the stationary solution of \eqref{intro:mainSDEsystem} induces a coupling of the two probability measures $m_1(x)\,\mathrm{d}x$ and $m_2(y)\,\mathrm{d}y$. The additional structural assumption $h \equiv f^2$ leads to a remarkable phenomenon:
The stationary solution of \eqref{intro:mainSDEsystem} is unique and, moreover, $X_t$ and $Y_t$ are independent for each $t \ge 0$. The joint density of $(X_t,Y_t)$ is thus given explicitly by the product $m_1(x)\,m_2(y)$ of the two marginal densities, regardless of the choice of the function $f$ (within the class permitted by Assumption \ref{asmp:main}). Note of course that $(X_t)_{t \ge 0}$ and $(Y_t)_{t \ge 0}$ cannot be independent \emph{as processes} for non-constant $f$, because $(Y_t)_{t \ge 0}$ appears in the dynamics of $(X_t)_{t \ge 0}$. The proof and additional discussion of this phenomenon are given in Section \ref{se:3}.

\begin{theorem}\label{thm:uniq}
Suppose that Assumption \ref{asmp:main} holds, and $h \equiv f^2$. Then the solution of the SDE \eqref{intro:mainSDEsystem} with the stationarity property $(X_t,Y_t)\stackrel{d}{=}(X_0,Y_0)$ for all $t\ge0$ is pathwise unique. That is, it admits a strong solution which is unique in law. In addition, the solution is independent in the sense that $X_t$ and $Y_t$ are independent for each $t \ge 0$. That is, for each $t \ge 0$, the law of $(X_t,Y_t)$ admits the density $m_1(x)\,m_2(y)$, where $m_1$ and $m_2$ are those of \eqref{def:marginals}.
\end{theorem}

The SDE \eqref{intro:mainSDEsystem} can be viewed as a McKean--Vlasov SDE (a.k.a.\ a non-linear SDE) in the sense that its coefficients depend not only on time and the current value of the solution but also on the fixed-time marginal distribution ${\mathcal L}(X_t,Y_t)$ of the solution at the time in consideration. The main challenge, in comparison with the classical theory of McKean--Vlasov SDEs (see, e.g., \cite{gartner1988mckean}, \cite{sznitman1991topics}, and the references in the latter) lies in the presence of the \textit{conditional} expectations in  \eqref{intro:mainSDEsystem}, most importantly in the diffusion coefficient of $(X_t)_{t\ge0}$. The underlying operation of passing from the joint distribution ${\mathcal L}(X_t,Y_t)$ to the conditional distribution ${\mathcal L}(Y_t|X_t)$ is notoriously discontinuous with respect to the weak convergence of probability measures. McKean--Vlasov SDEs involving conditional expectations have also been considered recently by \textsc{Jourdain}, \textsc{Leli\`{e}vre}, \textsc{Rousset}, \textsc{Roux} and \textsc{Stoltz} \cite{JLR}, \cite{lelievre2008long}, by \textsc{Bossy}, \textsc{Jabir} and \textsc{Talay} \cite{bossy2011conditional}, \cite{bossy2018wellposedness}, and by \textsc{Dermoune} \cite{dermoune}, who were respectively interested in the efficient simulation of Gibbs measures, turbulent flows, and adhesion particle dynamics. In contrast to our setup, in the SDEs of \cite{JLR}, \cite{lelievre2008long}, and \cite{dermoune}, a conditional expectation enters only into the drift coefficient of $(X_t)_{t\ge0}$, whereas in the SDEs of \cite{bossy2011conditional}, \cite{bossy2018wellposedness} the conditional distribution ${\mathcal L}(X_t|Y_t)$ rather than ${\mathcal L}(Y_t|X_t)$ enters into the coefficients of $(X_t)_{t\ge0}$. It is also worth mentioning that the weak uniqueness of the stationary solution can fail for McKean--Vlasov SDEs (see \cite{HT} and the references therein), which renders the weak uniqueness of the stationary solution for the SDE \eqref{intro:mainSDEsystem} in the full generality of Assumption \ref{asmp:main} an intriguing open problem.

\medskip

Despite the considerable interest in an inversion of the Markovian projection, very few rigorous results for the SDE \eqref{intro:mainSDEsystem} and its variants have been established so far. The paper by \textsc{Abergel} and \textsc{Tachet} \cite{abergel2010nonlinear} proves the existence in small time for forward Kolmogorov equations satisfied by the fixed-time marginal distributions ${\mathcal L}(X_t,Y_t)$ arising from SDEs like \eqref{intro:mainSDEsystem}, allowing for a multidimensional diffusion $(Y_t)_{t\ge0}$ and correlated $(W,B)$ but imposing a restrictive and somewhat implicit smallness assumption on $f'$. \textsc{Jourdain} and \textsc{Zhou} \cite{jourdain2016existence} showed the weak existence for a variant of the SDE \eqref{intro:mainSDEsystem} in which the diffusion $(Y_t)_{t\ge0}$ is replaced by a finite-state continuous-time Markov chain, assuming an insightful yet mysterious structural condition on the range of $f$ (\cite[Condition (C)]{jourdain2016existence}). Less closely related, \textsc{Alfonsi}, \textsc{Labart} and \textsc{Lelong} \cite{alfonsi2016stochastic} establish the strong existence and uniqueness of a conceptually similar counting process $(X_t)_{t\ge0}$ whose jump intensity involves a conditional expectation akin to the ones in \eqref{intro:mainSDEsystem}.  

\medskip

Solving the SDE \eqref{intro:mainSDEsystem} (or, more generally, the SDE \eqref{eq:SDE_gen}) allows one to construct processes that mimic the fixed-time marginal distributions of a given one-dimensional diffusion. When $b_1\equiv0$, this problem can be put into the broader context of martingale constructions with given fixed-time marginal distributions (a.k.a.\ peacocks). The latter have received much attention in stochastic analysis and financial mathematics. We refer the interested reader to the book \cite{hirsch2011peacocks} and the references therein. 

\subsection{Application to local stochastic volatility modeling}

A particularly prominent application of Proposition \ref{th:intro:Markovprojection} has been to the calibration of local stochastic volatility models (see \cite{lipton2002masterclass}, \cite{liptonmcghee2002masterclass}, \cite{piterbarg2006markovian}, \cite{guyon2011solved}, \cite{guyon2012being}, \cite[Chapter 11]{guyon2013nonlinear}, \cite{tian2015calibrating}, \cite{abergel2017nonparametric}, \cite{saporito2017calibration}). Hereby, one starts with a stochastic volatility process $(Z_t)_{t\ge0}$ (popular choices being the Cox--Ingersoll--Ross and exponential Ornstein--Uhlenbeck processes) and models the risk-neutral price $(S_t)_{t\ge0}$ of an asset by 
\begin{align}
\mathrm{d}S_t = S_tZ_t \sigma(t,S_t)\,\mathrm{d}W_t, \label{intro:SVmodel}
\end{align}
where $(W_t)_{t\ge0}$ is a Wiener process and $\sigma$ a function to be determined. Combinations of local and stochastic volatility models in this form have become 
quite popular and typically go by the name of \emph{local stochastic volatility (LSV)} \textit{models}. The stochastic process $(Z_t)_{t\ge0}$ provides greater flexibility than plain local volatility models, which are flexible enough to perfectly fit the implied volatility surface but require frequent re-calibration and typically fail to adequately incorporate exotic risks. We point out the implicit assumptions in  \eqref{intro:SVmodel} that the interest rate is zero and the asset pays no dividends.

\medskip

According to Proposition \ref{th:intro:Markovprojection}, if we define 
\begin{align}
\widehat{\sigma}(t,x) = \sigma(t,x)\sqrt{\E[Z_t^2|S_t=x]}, \label{intro:SVmodel-projection}
\end{align}
then the LSV model \eqref{intro:SVmodel} leads to the same fixed-time marginal distributions as the local volatility model
\begin{align}
\mathrm{d}\widehat{S}_t = \widehat{S}_t\widehat{\sigma}(t,\widehat{S}_t)\,\mathrm{d}\widehat{W}_t. \label{intro:LVmodel}
\end{align}
In particular, both models \eqref{intro:SVmodel} and \eqref{intro:LVmodel} produce the same European option prices. On the other hand, Dupire \cite{dupire1994pricing} famously showed that in order to be exactly calibrated to the observed call option prices $\{C(t,K):\,t > 0,\,K > 0\}$ the local volatility model \eqref{intro:LVmodel} must satisfy
\[
\widehat{\sigma}^2(t,K) = \frac{2 \partial_tC(t,K)}{K^2 \partial_{KK}C(t,K)} =: \sigma_{\mathrm{Dup}}^2(t,K).
\]
Hence, the original LSV model \eqref{intro:SVmodel} is exactly calibrated to the observed prices if
\[
\sigma(t,S_t) = \frac{\sigma_{\mathrm{Dup}}(t,S_t)}{\sqrt{\E[Z_t^2 | S_t]}}.
\]
Plugging this into \eqref{intro:SVmodel} yields the SDE
\begin{align}
\mathrm{d}S_t = S_t Z_t \frac{\sigma_{\mathrm{Dup}}(t,S_t)}{\sqrt{\E[Z_t^2 | S_t]}}\,\mathrm{d}W_t \label{intro:mainSLV-SDE}
\end{align}
or, equivalently, in terms of the log-price process $(X_t)_{t\ge0} = (\log S_t)_{t\ge0}$,
\begin{align}
\mathrm{d}X_t = -\frac12 Z_t^2 \frac{\sigma_{\mathrm{Dup}}^2(t,e^{X_t})}{\E[Z_t^2|X_t]}\,\mathrm{d}t + Z_t \frac{\sigma_{\mathrm{Dup}}(t,e^{X_t})}{\sqrt{\E[Z_t^2|X_t]}}\,\mathrm{d}W_t. \label{intro:mainSLV-SDElog}
\end{align}
Exact calibration of the LSV model thus reduces to constructing a solution of \eqref{intro:mainSLV-SDE} (or \eqref{intro:mainSLV-SDElog}), whose existence has been described as both ``a common belief in the quant community" \cite[p.\ 301]{guyon2013nonlinear} and a ``very challenging and open" problem \cite[p.\ 274]{guyon2013nonlinear}.

\medskip

One usually chooses $Z_t=f(Y_t)$, $t\ge0$, with an auxiliary one-dimensional diffusion $(Y_t)_{t\ge0}$ as in \eqref{eq:Y_def}. Moreover, only some of the call option prices $\{C(t,K):\,t>0,\,K > 0\}$ are available in reality, so that at the time of calibration one may choose to obtain a time-independent local volatility estimate $\widehat{\sigma}^2_{\mathrm{Dup}}(\cdot)$ instead of a time-dependent $\widehat{\sigma}^2_{\mathrm{Dup}}(t,\cdot)$. With these two choices, \eqref{intro:mainSLV-SDElog} falls exactly into the framework of \eqref{intro:mainSDEsystem}, notably with $h \equiv f^2$ as in Theorem \ref{thm:uniq}. For technical reasons, we need $W$ and $B$ to be uncorrelated and the coefficient $b_1$ in the SDE \eqref{intro:mainSDEsystem} to decrease (increase resp.) linearly as $x\to\infty$ ($x\to-\infty$ resp.), which requires a modification of $\widehat{\sigma}^2_{\mathrm{Dup}}(e^x)$ in the drift coefficient of \eqref{intro:mainSLV-SDElog} for large absolute values of the log-price $x$. Apart from these limitations, Theorems \ref{th:intro:main} and \ref{thm:uniq} are the first global existence and uniqueness results, respectively, in the context of the calibrated LSV model.  Results of this kind are of major importance, in particular, to ensure the accuracy of the widely used numerical solutions that have been proposed for the calibrated LSV model, such as the particle approximation method of \textsc{Guyon} and \textsc{Henry-Labord\`ere} \cite{guyon2012being}, as well as the (regularized) finite-difference approximation schemes with alternating directions in \cite{tian2015calibrating}, \cite{abergel2017nonparametric}, \cite{saporito2017calibration}. We refer to \cite[Chapter 11]{guyon2013nonlinear} for a detailed development of the LSV model and the associated calibration problem.

\subsection{Outline of the paper}

The rest of the paper is structured as follows. Section \ref{se:2.1} collects various important ingredients for the proofs of Theorems \ref{th:intro:main} and \ref{thm:uniq}. Specifically, in Subsection \ref{se:2_2} we describe a transformation of the stationary Fokker-Planck equation corresponding to the SDE \eqref{intro:mainSDEsystem}, introduced in Subsection \ref{se:2_1}, which enables us to apply the regularity estimate for invariant measures of \textsc{Bogachev}, \textsc{Krylov} and \textsc{R\"ockner} \cite[Theorem 1.1]{bogachev1996regularity} in Subsection \ref{se:2_4}. Subsection \ref{se:2_3} is devoted to the probabilistic analysis of the SDEs resulting from our transformation. In Section \ref{se:2.2} we deduce Theorem \ref{th:intro:main} from the results of Section \ref{se:2.1} by means of \textsc{Schauder}'s fixed-point theorem and \textsc{Veretennikov}'s pathwise uniqueness theorem for one-dimensional SDEs. Finally, Section \ref{se:3} provides the proof of Theorem \ref{thm:uniq}, which is based on another transformation of the stationary Fokker-Planck equation associated with the SDE \eqref{intro:mainSDEsystem}. 

\medskip

\noindent\textbf{Acknowledgements.} We thank Julien Guyon, Kasper Larsen, Alexander Lipton, and Scott Robertson for numerous enlightening discussions and references. We are especially grateful to Scott Robertson for introducing us to the subject of LSV models and for calling our attention to the interesting special case which became Theorem \ref{thm:uniq}.


\newcommand{\ep}{\hfill \ensuremath{\Box}}

\section{Ingredients for the proofs of Theorems \ref{th:intro:main} and \ref{thm:uniq}} \label{se:2.1}

\subsection{The stationary Fokker-Planck equation} \label{se:2_1}

The starting point for the proofs of Theorems \ref{th:intro:main} and \ref{thm:uniq} is the stationary Fokker-Planck equation associated with the SDE \eqref{intro:mainSDEsystem}. To state the latter we define, for any measurable function $\psi\!:\R\to(0,\infty)$ and probability density function $p$ on $\R^2$, the measurable function
\begin{equation*}
G^{\psi;p}:\;\;\R\to[0,\infty),\quad x\mapsto \frac{\int_{\R}p(x,y)\,\mathrm{d}y}{\int_{\R}\psi(y)\,p(x,y)\,\mathrm{d}y}.
\end{equation*}
Notice that, if $(U,V)$ is a random vector with the joint probability density function $p$, then $G^{\psi;p}(x) = \E[\psi(V)|U\!=\!x]^{-1}$. Putting this together with Dynkin's formula we conclude that any probability density function $p$ stationary for the SDE \eqref{intro:mainSDEsystem} must be a generalized solution of the stationary Fokker-Planck equation   
\begin{equation}\label{mainpde}
\begin{split}
 0=\frac12\big(\sigma_1^2(x)\,f^2(y)\,G^{f^2;p}(x)\,p(x,y)\big)_{xx} + \frac12\big(\sigma_2^2(y)\,p(x,y)\big)_{yy} \\
 -\big(b_1(x)\,h(y)\,G^{h;p}(x)\,p(x,y)\big)_x-\big(b_2(y)\,p(x,y)\big)_y.
\end{split}
\end{equation}
More specifically, each stationary $p$ has the property that
\begin{equation}\label{mainpde_weak}
\begin{split}
0 = \int_{\R^2}\bigg(&\frac12\,\sigma_1^2(x)\, f^2(y)\,G^{f^2;p}(x)\,\varphi_{xx}(x,y) + \frac12\, \sigma_2^2(y)\,\varphi_{yy}(x,y) \qquad \quad\;\; \\ 
&+\,b_1(x)\,h(y)\,G^{h;p}(x)\,\varphi_x(x,y) + b_2(y)\,\varphi_y(x,y) 
\!\bigg) p(x,y)\,\mathrm{d}x\,\mathrm{d}y, \quad \varphi \in C^\infty_c(\R^2),
\end{split}
\end{equation}
where $C^\infty_c(\R^2)$ is the space of all infinitely differentiable functions from $\R^2$ to $\R$ with compact support.

\subsection{The main transformation} \label{se:2_2}

The conditional expectation terms $G^{f^2;p}$, $G^{h;p}$ render the PDE \eqref{mainpde} non-linear and involve both local and non-local effects. However, if we think of $G^{f^2;p}$, $G^{h;p}$ as given, the PDE \eqref{mainpde} becomes a linear stationary Fokker-Planck equation, as studied extensively in \cite{bogachev1996regularity}. We recall their main result in the finite-dimensional case, which allows to control the gradient of a solution $p$ in terms of some but, crucially, not all first-order derivatives of the diffusion matrix. Hereby, with $d\ge1$, we use the notation $C_b^\infty(\R^d)$ for the space of infinitely differentiable functions from $\R^d$ to $\R$ with bounded derivatives of all orders and $W^1_2(\R^d)$ for the Sobolev space of functions from $\R^d$ to $\R$ that are square integrable together with their gradients. 

\begin{proposition}[cf.\ \cite{bogachev1996regularity}, Theorem 1.1]\label{BKR_reg}
Let $\alpha>0$, $A=(A^{ij})_{1\le i,j\le d}$ be a bounded uniformly Lipschitz continuous function from $\R^d$ to the set of the symmetric $d\times d$ matrices whose eigenvalues are bounded below by $\alpha$, and $B=(B^i)_{1\le i\le d}$ be a measurable function from $\R^d$ to $\R^d$. With 
\begin{equation*}
{\mathcal L} =\sum_{i,j=1}^d A^{ij}\,\partial_{z_iz_j}+\sum_{i=1}^d B^i\,\partial_{z_i},
\end{equation*}
consider a probability measure $\mu$ on $\R^d$ such that $\int_{\R^d} |B|^2\,\mathrm{d}\mu < \infty$ and 
\begin{equation*}
\forall \varphi \in C^\infty_b(\R^d):\quad \int_{\R^d} {\mathcal L}\varphi\,\mathrm{d}\mu=0.
\end{equation*}
Then $\mu(\mathrm{d}z)=p(z)\,\mathrm{d}z$, with $\sqrt{p}\in W^1_2(\R^d)$, and 
\begin{equation*}
4\int_{\R^d} \big|\nabla\sqrt{p}\big|^2\,\mathrm{d}z
=\int_{\R^d} \left|\frac{\nabla p}{p}\right|^2\,p\,\mathrm{d}z
\leq \dfrac{2}{\alpha^2}\int_{\R^d}\big(|B|^2+|D|^2\big)\,p\,\mathrm{d}z,
\end{equation*}
where $D=\big(\sum_{i=1}^d \partial_{z_i} A^{ij}\big)_{1\le j\le d}$ and we have adopted the convention that $\frac{\nabla p}{p}\equiv 0$ outside of the support of $p$. 
\end{proposition}

In the case of \eqref{mainpde_weak}, the diffusion matrix $A$ is given by 
\[
\begin{pmatrix}
\frac12\,\sigma_1^2(x)\,f^2(y)\,G^{f^2;p}(x) & 0 \\
0 & \frac12\,\sigma_2^2(y)
\end{pmatrix}
= \begin{pmatrix}
\frac12\,\sigma_1^2(x)\,\frac{f^{2}(y)}{\E[f^{2}(V)|U=x]} & 0 \\
0 & \frac12\,\sigma_2^2(y)
\end{pmatrix},
\]
where $(U,V)$ is a random vector with the joint probability density function $p$. Thus, a direct use of Proposition \ref{BKR_reg} would require an a priori regularity estimate on the function $x\mapsto\E[f^{2}(V)|U\!=\!x]$. We circumvent this difficulty by applying a suitable transformation to $p$. The latter acts on the space $L^1_{\mathrm{prob}}(\R^2)$ of probability density functions on $\R^2$ and is defined by 
\begin{equation}\label{meas_change}
T:\;\;L^1_{\mathrm{prob}}(\R^2)\to L^1_{\mathrm{prob}}(\R^2),\quad p\mapsto
p(x,y)\,f^2(y)\,G^{f^2;p}(x).
\end{equation} 
The key observation, established in the next proposition, is that, if \eqref{mainpde_weak} holds for $p$, then $\widetilde{p}=Tp$ satisfies 
\begin{equation}\label{mainpde_weak-transformed}
\begin{split}
0 \!=\! \int_{\R^2}\! \bigg(\frac12\,\sigma_1^2(x)\,\varphi_{xx}(x,y) + \frac12\,\sigma_2^2(y)\, f^{-2}(y)\,G^{f^{-2};\widetilde{p}}(x)\,\varphi_{yy}(x,y)
\qquad\qquad\qquad\qquad\qquad\qquad\; \\
+ b_1(x) (hf^{-2})(y)\,G^{hf^{-2};\widetilde{p}}(x)\,\varphi_x(x,y)
\!+\! b_2(y)\,f^{-2}(y)\,G^{f^{-2};\widetilde{p}}(x)\,\varphi_y(x,y) 
\!\bigg) \widetilde{p}(x,y)\,\mathrm{d}x\,\mathrm{d}y
\end{split}
\end{equation}
for all $\varphi \in C^\infty_c(\R^2)$. At the first glance the equation \eqref{mainpde_weak-transformed} does not look any better than the equation \eqref{mainpde_weak} but they differ in one critical way: In \eqref{mainpde_weak-transformed}, the diffusion matrix reads
\[
\begin{pmatrix}
\frac12\,\sigma_1^2(x) & 0 \\
0 & \frac12\,\sigma_2^2(y)\,f^{-2}(y)\,G^{f^{-2};\widetilde{p}}(x)
\end{pmatrix}
=\begin{pmatrix}
\frac12\,\sigma_1^2(x) & 0 \\
0 & \frac12\,\sigma_2^2(y)\,\frac{f^{-2}(y)}{\E[f^{-2}(\widetilde{V})|\widetilde{U}=x]} 
\end{pmatrix},
\]
where $(\widetilde{U},\widetilde{V})$ is a random vector with the joint probability density function $\widetilde{p}$. Thus, in contrast to \eqref{mainpde_weak}, the derivatives $\partial_x A^{11}$, $\partial_y A^{22}$ do not involve the derivative of a conditional expectation term.   

\medskip

The following proposition summarizes the main properties of the transformation $T$. 

\begin{proposition}\label{prop_changemeas}
Under Assumption \ref{asmp:main} the following are true:	
\begin{enumerate}
\item[\rm{(i)}] For each $p\in L^1_{\mathrm{prob}}(\R^2)$, one has $\widetilde{p}=Tp\in L^1_{\mathrm{prob}}(\R^2)$ and the first marginal of $\widetilde{p}$ is the same as that of $p$, i.e.,\ $\int_{\R} p(\cdot,y)\,\mathrm{d}y = \int_{\R} \widetilde{p}(\cdot,y)\,\mathrm{d}y$.  
\item[\rm{(ii)}] $T$ is a bijection from $L^1_{\mathrm{prob}}(\R^2)$ to itself, and we have 
\begin{equation}\label{meas_change_back}
p(x,y)=(Tp)(x,y)\,f^{-2}(y)\,G^{f^{-2};Tp}(x),\quad p\in L^1_{\mathrm{prob}}(\R^2).
\end{equation}
\item[\rm{(iii)}] For any measurable function $\psi\!:\R\to(0,\infty)$ and $p\in L^1_{\mathrm{prob}}(\R^2)$, it holds 
\[
G^{\psi;p}=\frac{G^{\psi f^{-2};Tp}}{G^{f^{-2};Tp}}.
\]
\item[\rm{(iv)}]\label{c} $p$ satisfies \eqref{mainpde_weak} if and only if $\widetilde{p}=Tp$ satisfies \eqref{mainpde_weak-transformed} for all $\varphi \in C^\infty_c(\R^2)$.
\end{enumerate}
\end{proposition}

\smallskip

\noindent\textbf{Proof.\ (i).} Inserting the definition of $G^{f^2;p}$ into $\widetilde{p}(x,y)=p(x,y)\,f^2(y)\,G^{f^2;p}(x)$ and integrating in $y$ we obtain
\begin{equation*}
\int_\mathbb{R} \widetilde{p}(x,y)\,\mathrm{d}y
=\int_\mathbb{R} f^2(y)\,p(x,y)\,\mathrm{d}y\,\dfrac{\int_\mathbb{R} p(x,y)\,\mathrm{d}y}{\int_\mathbb{R} f^2(y)\,p(x,y)\mathrm{d}y}=\int_\mathbb{R} p(x,y)\,\mathrm{d}y.
\end{equation*}
In particular, $\int_{\R} \int_{\R} \widetilde{p}(x,y)\,\mathrm{d}y\,\mathrm{d}x=1$.

\medskip

\noindent\textbf{(ii).} Plugging the definition of $G^{f^2;p}$ into $f^{-2}(y)\,\widetilde{p}(x,y)=p(x,y)\,G^{f^2;p}(x)$ and integrating in $y$ we get
\[
\int_{\mathbb{R}} f^{-2}(y)\,\widetilde{p}(x,y)\,\mathrm{d}y
=\frac{\big(\int_{\mathbb{R}} p(x,y)\,\mathrm{d}y\big)^2}{\int_{\mathbb{R}} f^2(y)\,p(x,y)\,\mathrm{d}y}.
\]
Combining the latter equation with part (i) we arrive at
\[
\frac{\int_{\R} f^{-2}(y)\,\widetilde{p}(x,y)\,\mathrm{d}y}
{\int_{\R} \widetilde{p}(x,y)\,\mathrm{d}y}
=\frac{\int_{\R} p(x,y)\,\mathrm{d}y}{\int_{\R} f^2(y)\,p(x,y)\,\mathrm{d}y},
\]
that is, $G_{f^2;p}=1/G^{f^{-2};\widetilde{p}}$. The identity \eqref{meas_change_back} and, in particular, the injectivity of $T$ then follow by rearranging $\widetilde{p}(x,y)=p(x,y)\,f^2(y)\,/\,G^{f^{-2};\widetilde{p}}(x)$. Moreover, for any $\widetilde{p}\in L^1_{\mathrm{prob}}(\R^2)$, the probability density function $p(x,y)=\widetilde{p}(x,y)\,f^{-2}(y)\,G^{f^{-2};\widetilde{p}}(x)$ satisfies 
\begin{align*}
(Tp)(x,y) &=\widetilde{p}(x,y)\,G^{f^{-2};\widetilde{p}}(x)\,G^{f^2;p}(x)
=\widetilde{p}(x,y)\,G^{f^{-2};\widetilde{p}}(x)\,\frac{\int_{\R} \widetilde{p}(x,y)\,f^{-2}(y)\,G^{f^{-2};\widetilde{p}}(x)\,\mathrm{d}y}{\int_{\R} \widetilde{p}(x,y)\,G^{f^{-2};\widetilde{p}}(x)\,\mathrm{d}y} \\
&=\widetilde{p}(x,y),
\end{align*}
which shows the surjectivity of $T$. 

\medskip

\noindent\textbf{(iii).} Writing $\widetilde{p}$ for $Tp$ as before. Using the definition of $G^{\psi;p}$ as well as parts (i) and (ii), we compute 
\[
G^{\psi;p}(x) = \frac{\int_{\R} p(x,y)\,\mathrm{d}y}{\int_{\R} \psi(y)\,p(x,y)\,\mathrm{d}y} 
= \frac{\int_{\R} \widetilde{p}(x,y)\,\mathrm{d}y}{\int_{\R}\psi(y)\,\widetilde{p}(x,y)\,f^{-2}(y)\,G^{f^{-2};\widetilde{p}}(x)\,\mathrm{d}y} 
= \frac{G^{\psi f^{-2};\widetilde{p}}(x)}{G^{f^{-2};\widetilde{p}}(x)}.
\]

\smallskip

\noindent\textbf{(iv).} The definition of $\widetilde{p}=Tp$ reveals that the equation in \eqref{mainpde_weak} is equivalent to 
\[
\begin{split}
0 = \int_{\R^2} \bigg(\frac12\,\sigma_1^2(x)\,\varphi_{xx}(x,y) + \frac{\sigma_2^2(y)}{2f^2(y)\,G^{f^2;p}(x)}\,\varphi_{yy}(x,y) 
\qquad\qquad\qquad\qquad\qquad\;\; \\
+\frac{b_1(x)\,h(y)\,G^{h;p}(x)}{f^2(y)\,G^{f^2;p}(x)}\,\varphi_x(x,y)
 + \frac{b_2(y)}{f^2(y)\,G^{f^2;p}(x)}\,\varphi_y(x,y)
\!\bigg)\,\widetilde{p}(x,y)\,\mathrm{d}x\,\mathrm{d}y,
\end{split}
\]
for all $\varphi \in C^\infty_c(\R^2)$.
As observed in the proof of part (ii), it holds $1/G^{f^2;p}=G^{f^{-2};\widetilde{p}}$. Inserting also $G^{h;p}/G^{f^2;p}=G^{h;p}\,G^{f^{-2};\widetilde{p}}=G^{hf^{-2};\widetilde{p}}$ (cf.\ part (iii)) we end up with \eqref{mainpde_weak-transformed}. \ep

\subsection{The transformed SDE} \label{se:2_3}

We proceed by noting that any $\widetilde{p}\in L^1_{\mathrm{prob}}(\R^2)$ stationary for the SDE
\begin{align}
\begin{cases}
\;\mathrm{d}\widetilde{X}_t = b_1(\widetilde{X}_t)\,\frac{h(\widetilde{Y}_t)f^{-2}(\widetilde{Y}_t)}{\E[h(\widetilde{Y}_t) f^{-2}(\widetilde{Y}_t)|\widetilde{X}_t]}\,\mathrm{d}t + \sigma_1(\widetilde{X}_t)\,\mathrm{d}W_t, \vspace{3pt} \\
\;\mathrm{d}\widetilde{Y}_t \;= b_2(\widetilde{Y}_t)\,\frac{f^{-2}(\widetilde{Y}_t)}{\E[f^{-2}(\widetilde{Y}_t)| \widetilde{X}_t]}\,\mathrm{d}t + \sigma_2(\widetilde{Y}_t)\,\frac{f^{-1}(\widetilde{Y}_t)}{\sqrt{\E[f^{-2}(\widetilde{Y}_t)|\widetilde{X}_t]}}\,\mathrm{d}B_t
\end{cases} \label{intro:transformedSDEsystem}
\end{align}
must satisfy \eqref{mainpde_weak-transformed} for all $\varphi \in C^\infty_c(\R^2)$. The SDE \eqref{intro:transformedSDEsystem}, in turn, falls into the more general framework of the next proposition.  

\begin{proposition}\label{pr:general linear SDE}
Suppose the measurable functions $\overline{b}_1,\overline{b}_2,\overline{\sigma}_1,\overline{\sigma}_2:\,\R^2\to\R$ obey
\begin{equation}\label{cond_coeff}
\begin{split}
\forall\,x,y \in \R:\qquad & x\overline{b}_1(x,y) \le -\overline{c}x^2 + \overline{C}_1, \qquad\; y\overline{b}_2(x,y) \le -\overline{c}y^2 + \overline{C}_1, \\ 
& |\overline{b}_1(x,y)| \le \overline{C}_2(1+|x|), \qquad
|\overline{b}_2(x,y)| \le \overline{C}_2(1+|y|), \\
& 2\alpha \le \overline{\sigma}^2_1(x,y),\overline{\sigma}^2_2(x,y)\leq 2\Sigma, 
\end{split}
\end{equation}
with some constants $\overline{c},\overline{C}_1,\overline{C}_2,\alpha,\Sigma\in(0,\infty)$. Then the weak solution of the SDE
\begin{equation}
\begin{cases}
\;\mathrm{d}\overline{X}_t=\overline{b}_1(\overline{X}_t,\overline{Y}_t)\,\mathrm{d}t+\overline{\sigma}_1(\overline{X}_t,\overline{Y}_t)\,\mathrm{d}W_t,
 \\
\;\mathrm{d}\overline{Y}_t\,=\overline{b}_2(\overline{X}_t,\overline{Y}_t)\,\mathrm{d}t+\overline{\sigma}_2(\overline{X}_t,\overline{Y}_t)\,\mathrm{d}B_t
\end{cases} \label{lem:linear-SDE}
\end{equation} 
with the property $(\overline{X}_t,\overline{Y}_t)\stackrel{d}{=}(\overline{X}_0,\overline{Y}_0)$ for $t\ge0$ is unique in law. Moreover, it satisfies the following:
\begin{enumerate}
\item[\rm{(i)}] $\E\big[\overline{X}_t^2 + \overline{Y}_t^2\big] \leq (2\Sigma+\overline{C}_1)/\overline{c}$, for $t\ge0$.
\item[\rm{(ii)}] $\overline{X}_0$ has a density $m$ such that for each $R\in(0,\infty)$ there exists some $\delta_R>0$ with $m\ge\delta_R$ a.e.\ in $[-R,R]$, where $\delta_R$ can be chosen to depend on $R,\overline{c},\overline{C}_1,\overline{C}_2,\alpha,\Sigma$ only. 
\end{enumerate}
\end{proposition}

The proof makes use of a well-known result from \cite{khasminskii2011stochastic}. 

\begin{proposition}[cf.\ \cite{khasminskii2011stochastic}, Theorem 4.1 and Corollary 4.4] \label{pr:khasminskii}
Let $b:\,\R^d\to\R^d$ and $\sigma:\,\R^d\to\R^{d\times\ell}$ be measurable functions and $(Z_t)_{t\ge0}$ be a non-explosive time-homogeneous Markov process in $\R^d$ described by the SDE
\begin{equation*}
\mathrm{d}Z_t=b(Z_t)\,\mathrm{d}t + \sigma(Z_t)\,\mathrm{d}\beta_t,
\end{equation*}
where $(\beta_t)_{t\ge0}$ is a standard Brownian motion in $\R^\ell$. Define $A(z)=\sigma(z)\sigma(z)^\top$, $z\in\R^d$ and suppose that there exists a bounded open $U\subset \mathbb{R}^d$ with $C^1$--boundary having the properties:
\begin{enumerate}
\item[\rm{(a)}] The smallest eigenvalue of the diffusion matrix $A(z)$ is uniformly bounded away from zero on an open neighborhood of $U$.  
\item[\rm{(b)}] For each compact $K\subset\R^d$, it holds $\sup_{z\in K}\E[\tau|Z_0=z]<\infty$, where $\tau$ is the hitting time of the set $U$.
\end{enumerate}
Then the Markov process $(Z_t)_{t\ge0}$ admits a unique stationary distribution.
\end{proposition}

\smallskip

\noindent\textbf{Proof of Proposition \ref{pr:general linear SDE}.} We start by pointing out that the condition \eqref{cond_coeff} on the coefficients of the SDE \eqref{lem:linear-SDE} suffices to ensure that its associated martingale problem is well-posed. This follows from \cite[Theorem 10.2.2 and Exercise 7.3.4]{stroock-varadhan} (see also \cite[Theorem 3 and the paragraph following it]{krylov1969ito}). Consequently, there is a unique non-explosive time-homogeneous strong Markov process in $\R^2$ described by the SDE \eqref{lem:linear-SDE} (see, e.g., \cite[Chapter 5, Theorem 4.20]{karatzas-shreve}).

\medskip

Property (a) of Proposition \ref{pr:khasminskii} holds for the SDE \eqref{lem:linear-SDE} with any choice of $U$ by assumption, so we turn to checking property (b). Suppose $(\overline{X}_t,\overline{Y}_t)_{t\ge0}$ solves the SDE \eqref{lem:linear-SDE} for some fixed initial position. By It\^o's formula,
\begin{align}
\begin{split}
\mathrm{d}\overline{X}_t^2 &=\big(\overline{\sigma}_1^2(\overline{X}_t,\overline{Y}_t)+2\overline{X}_t\overline{b}_1(\overline{X}_t,\overline{Y}_t)\big)\,\mathrm{d}t+2\overline{X}_t\overline{\sigma}_1(\overline{X}_t,\overline{Y}_t)\,\mathrm{d}W_t, \\
\mathrm{d}\overline{Y}_t^2 &=\big(\overline{\sigma}_2^2(\overline{X}_t,\overline{Y}_t)+2\overline{Y}_t\overline{b}_2(\overline{X}_t,\overline{Y}_t)\big)\,\mathrm{d}t+2\overline{Y}_t\overline{\sigma}_2(\overline{X}_t,\overline{Y}_t)\,\mathrm{d}B_t.
\end{split} \label{pf:generallinear1}
\end{align}
Next, we let $R_t = \big(1+\overline{X}_t^2+\overline{Y}_t^2\big)^{1/2}$ and compute
\begin{align*}
& \mathrm{d}R_t =\frac{\overline{X}_t}{R_t}\,\overline{\sigma}_1(\overline{X}_t,\overline{Y}_t)\,\mathrm{d}W_t + \frac{\overline{Y}_t}{R_t}\,\overline{\sigma}_2(\overline{X}_t,\overline{Y}_t)\,\mathrm{d}B_t \\ &\qquad\;\; +\!\frac{1}{R_t}\bigg(\!\overline{\sigma}_1^2(\overline{X}_t,\overline{Y_t})\bigg(\frac12\!-\!\frac{\overline{X}_t^2}{R_t^2}\bigg) \!+\! \overline{\sigma}_2^2(\overline{X}_t,\overline{Y_t})\bigg(\frac12\!-\!
\frac{\overline{Y}_t^2}{R_t^2}\bigg) \!+\! \overline{X}_t\overline{b}_1(\overline{X}_t,\overline{Y}_t) \!+\! \overline{Y}_t\overline{b}_2(\overline{X}_t,\overline{Y}_t)\!\bigg) \mathrm{d}t.
\end{align*}
In the latter expression, the diffusion coefficients are bounded, whereas the drift coefficient is less than or equal to
\begin{align*}
\frac{1}{R_t}\big(2\Sigma - \overline{c}\overline{X}_t^2 - \overline{c}\overline{Y}_t^2 + 2\overline{C}_1\big) = \frac{2\Sigma+\overline{c}+2\overline{C}_1}{R_t} - \overline{c}R_t
\end{align*}
and, thus, negative and uniformly bounded away from zero whenever $R_t^2\ge\frac{2\Sigma+\overline{c}+2\overline{C}_1}{\overline{c}}+1$. With $U=\{(x,y)\in\R^2:\,x^2+y^2<(2\Sigma+\overline{c}+2\overline{C}_1)/\overline{c}\}$, we can now use a simple time-change argument relying on the Dambis-Dubins-Schwarz theorem (see, e.g., \cite[Chapter 3, Theorem 4.6]{karatzas-shreve} and recall the assumptions imposed on $\overline{\sigma}^2_1$, $\overline{\sigma}^2_2$ in \eqref{cond_coeff}) for the martingale part of $(R_t)_{t\ge0}$ to obtain property (b) of Proposition \ref{pr:khasminskii}.

\medskip

To prove claim (i) we write $(x_0,y_0)$ for the initial position of $(\overline{X}_t,\overline{Y}_t)_{t\ge0}$ and localize by means of the stopping times $\tau_n=\inf\{t\ge0:\, \overline{X}_t^2+\overline{Y}_t^2 \ge n\}$, $n\in\N$, for the case that the local martingale part is not a true martingale, deducing from \eqref{pf:generallinear1}: 
\begin{align*}
& \;\E\big[\overline{X}_{t \wedge \tau_n}^2 + \overline{Y}_{t \wedge \tau_n}^2\big] \\
& = x_0^2+y_0^2+\E\bigg[\int_0^{t \wedge \tau_n} \overline{\sigma}_1^2(\overline{X}_s,\overline{Y_s})+\overline{\sigma}_2^2(\overline{X}_s,\overline{Y_s}) + 2\overline{X}_s \overline{b}_1(\overline{X}_s,\overline{Y}_s)+ 2\overline{Y}_s\overline{b}_2(\overline{X}_s,\overline{Y}_s)\,\mathrm{d}s\bigg] \\
&\le x_0^2+y_0^2+2\E\bigg[\int_0^{t \wedge \tau_n} 2\Sigma  - \overline{c}( \overline{X}_s^2 +  \overline{Y}_s^2) + \overline{C}_1\,\mathrm{d}s\bigg]. 
\end{align*}
By Fatou's lemma and the monotone convergence theorem,
\begin{align*}
\E\big[\overline{X}_t^2 + \overline{Y}_t^2\big] \le x_0^2 + y_0^2 + 2\E\bigg[\int_0^t 2\Sigma + \overline{C}_1  - \overline{c}(\overline{X}_s^2 + \overline{Y}_s^2)\,\mathrm{d}s\bigg].
\end{align*}
Dividing by $t$ and taking the limit inferior on both sides we get, in particular,
\begin{align*}
0 & \le \liminf_{t\to\infty}\,\frac{1}{t}\,\E\bigg[\int_0^t 2\Sigma + \overline{C}_1  - \overline{c}(\overline{X}_s^2 + \overline{Y}_s^2)\,\mathrm{d}s\bigg] \\
& \le \liminf_{t\to\infty}\,\frac{1}{t}\,\E\bigg[\int_0^t \big(2\Sigma + \overline{C}_1  - \overline{c}(\overline{X}_s^2 + \overline{Y}_s^2)\big)\vee(-M)\,\mathrm{d}s\bigg],\quad M\in\N. 
\end{align*}
Thanks to \cite[Corollary 4.3]{khasminskii2011stochastic} the latter limit inferior can be evaluated to an integral with respect to the invariant distribution of $(\overline{X}_t,\overline{Y}_t)_{t\ge0}$, which yields claim (i) after passing to the limit $M\to\infty$ via the monotone convergence theorem. 

\medskip

For the proof of claim (ii) we let $(\overline{X}_t,\overline{Y}_t)_{t\ge0}$ be the stationary solution of the SDE \eqref{lem:linear-SDE} and apply Proposition \ref{th:intro:Markovprojection} to conclude that $\overline{X}_t \stackrel{d}{=}\widehat{X}_t$ for all $t\ge0$, where $(\widehat{X}_t)_{t\ge0}$ is a stationary solution of the SDE
\begin{align*}
& \mathrm{d}\widehat{X}_t = \widehat{b}(\widehat{X}_t)\,\mathrm{d}t + \widehat{\sigma}(\widehat{X}_t)\,\mathrm{d}W_t, \\
& \widehat{b}(x) = \E[\overline{b}_1(\overline{X}_t,\overline{Y}_t)|\overline{X}_t\!=\!x],\quad\text{and}
\quad 
\widehat{\sigma}^2(x) = \E[\overline{\sigma}^2_1(\overline{X}_t,\overline{Y}_t)| \overline{X}_t\!=\!x].
\end{align*}
Then, for all $x\in\R$, we have $2\alpha \le \widehat{\sigma}^2(x) \le 2\Sigma$, as well as
\begin{align*}
& x\widehat{b}(x) = \E[\overline{X}_t\overline{b}_1(\overline{X}_t,\overline{Y}_t) | \overline{X}_t\!=\!x] \le \E[-\overline{c}\overline{X}_t^2 + \overline{C}_1 | \overline{X}_t\!=\!x] = -\overline{c}x^2+\overline{C}_1, \\
& |\widehat{b}(x)| \le \E[|\overline{b}_1(\overline{X}_t,\overline{Y}_t)|\,|\, \overline{X}_t\!=\!x] \le \E[\overline{C}_2(1+|\overline{X}_t|)|\overline{X}_t\!=\!x] 
= \overline{C}_2(1 + |x|).
\end{align*}
It follows that the (common) law $m(\mathrm{d}x)$ of $\overline{X}_t$, $t\ge0$ satisfies
\begin{align*}
\forall \varphi \in C^\infty_c(\R):\quad
\int_\R \frac12\,\widehat{\sigma}^2(x)\,\varphi''(x)+\widehat{b}(x)\,\varphi'(x) \,m(\mathrm{d}x)=0.
\end{align*}
In other words, $\tfrac12(\widehat{\sigma}^2 m)'' -(\widehat{b}m)' = 0$ in the sense of distributions. Consequently, it holds $\big(\frac12\widehat{\sigma}^2 m-\int_0^\cdot \widehat{b}\,\mathrm{d}m\big)''=0$ in the sense of distributions and, by elliptic regularity (see, e.g., \cite[Appendix B, Theorem 14]{lax}), in the classical sense. Therefore, $\frac12\widehat{\sigma}^2(x)\,m(\mathrm{d}x)-\int_0^x \widehat{b}(a)\,m(\mathrm{d}a)=k_1x+k_2$ for some $k_1,k_2\in\R$. This identity shows that $\frac12\widehat{\sigma}^2(x)\,m(\mathrm{d}x)$ is given by a locally bounded measurable function a priori and by a locally Lipschitz function a posteriori, which we denote by $\theta$. By differentiating we obtain
\begin{equation*}
\theta'(x)-\frac{2\widehat{b}(x)}{\widehat{\sigma}^2(x)}\,\theta(x)=k_1\quad\text{for a.e. }x\in\R. 
\end{equation*} 
Multiplying both sides by the integrating factor $e^{-\int_0^x \frac{2\widehat{b}(a)}{\widehat{\sigma}^2(a)}\,\mathrm{d}a}$, integrating the resulting equation, and rearranging we arrive at
\begin{equation*}
\theta(x)=\theta(0)\exp\bigg(\int_0^x \frac{2\widehat{b}(a)}{\widehat{\sigma}^2(a)}\,\mathrm{d}a\bigg)
+k_1\exp\bigg(\int_0^x \frac{2\widehat{b}(a)}{\widehat{\sigma}^2(a)}\,\mathrm{d}a\bigg)
\int_0^x \exp\bigg(-\int_0^{a_1} \frac{2\widehat{b}(a_2)}{\widehat{\sigma}^2(a_2)}\,\mathrm{d}a_2\bigg)\,\mathrm{d}a_1. 
\end{equation*}
The right-hand side defines the density of a finite positive measure only if $k_1=0$ and $\theta(0)>0$. Claim (ii) readily follows. \ep 

\subsection{The transformed PDE} \label{se:2_4}

This subsection is devoted to the analysis of the  stationary Fokker-Planck equation \eqref{mainpde_weak-transformed} in which the probability density function $\widetilde{p}$ within the non-linear non-local terms $G^{f^{-2};\widetilde{p}}$, $G^{hf^{-2};\widetilde{p}}$ is thought of as given.   

\begin{proposition} \label{pr:linearized-PDE-estimates}
Let $q \in L^1_{\mathrm{prob}}(\R^2)$. Under Assumption \ref{asmp:main} there exists a unique probability measure $\mu$ on $\R^2$ with a finite first moment such that 
\begin{align} \label{PDE-linearized}
\begin{split}
0 \!=\!\int_{\R^2} \! \bigg(\frac12\,\sigma_1^2(x)\,\varphi_{xx}(x,y) + \frac12\,\sigma_2^2(y)\, f^{-2}(y)\,G^{f^{-2};q}(x)\,\varphi_{yy}(x,y) 
\qquad\qquad \qquad\qquad\qquad\quad\;\;\, \\
+\,b_1(x)\,(hf^{-2})(y)\,G^{hf^{-2};q}(x)\,\varphi_x(x,y) 
+ b_2(y)\,f^{-2}(y)\,G^{f^{-2};q}(x)\,\varphi_y(x,y)\!\bigg) \mu(\mathrm{d}x,\mathrm{d}y), \\ 
\varphi \in C^\infty_c(\R^2).
\end{split}
\end{align}
Further, there are constants $\widetilde{C}_1,\widetilde{C}_2<\infty$ and, for any fixed $R\in(0,\infty)$, a constant $\widetilde{\delta}_R>0$, all depending only on the constants mentioned in Assumption \ref{asmp:main} (and, in particular, independent of $q$), such that 
\begin{enumerate}
\item[\rm{(i)}] $\mu(\mathrm{d}x,\mathrm{d}y) = \widetilde{q}(x,y)\,\mathrm{d}x\,\mathrm{d}y$ for some $\widetilde{q} \in L^1_{\mathrm{prob}}(\R^2)$ satisfying $\sqrt{\widetilde{q}} \in W^1_2(\R^2)$ and
\[
\frac14\int_{\R^2}\frac{|\nabla\widetilde{q}(x,y) |^2}{\widetilde{q}(x,y)}\,\mathrm{d}x\,\mathrm{d}y = \int_{\R^2}\big|\nabla \sqrt{\widetilde{q}}(x,y)\big|^2\,\mathrm{d}x\,\mathrm{d}y \le \widetilde{C}_1.
\]
\item[\rm{(ii)}] $\int_{\R^2}(x^2+y^2)\,\mu(\mathrm{d}x,\mathrm{d}y) \le \widetilde{C}_2$. 
\vspace{6pt}
\item[\rm{(iii)}] $\widetilde{m}(x) = \int_{\R}\widetilde{q}(x,y)\,\mathrm{d}y\ge\delta_R$ a.e.\ in $[-R,R]$. 
\end{enumerate}
\end{proposition}

An immediate application of Proposition \ref{BKR_reg} to the equation in \eqref{PDE-linearized} is hindered by the lack of an a priori regularity estimate on $G^{f^{-2};q}$. To address this, we mollify as follows. We pick a non-negative $\kappa\in C^\infty_c(\R)$ supported in $[-1,1]$ and such that $\int_\R \kappa(x)\,\mathrm{d}x=1$. With $\kappa_n(x)=n\kappa(nx)$, we set
\begin{align*}
G_n^{\psi;q}(x)=\int_{\R} \kappa_n(x-a)\,G^{\psi;q}(a)\,\mathrm{d}a
\end{align*}
for any measurable function $\psi:\,\R\to(0,\infty)$, $q\in L^1_{\mathrm{prob}}(\R^2)$ and $n\in\N$. If $c_\psi\le\psi\le C_\psi$ for some  $c_\psi,C_\psi\in(0,\infty)$, then $C_\psi^{-1}\le G^{\psi;q}\le c_\psi^{-1}$ and $C_\psi^{-1}\le G^{\psi;q}_n\le c_\psi^{-1}$	 for $n\in\N$. Also,
\begin{align}
\forall n\in\N:\quad |(G_n^{\psi;q})'| \le nc_\psi^{-1} \int_\R |\kappa'(x)|\,\mathrm{d}x. \label{def:mollified-gradbound}
\end{align}
The next lemma deals with a mollified version of the equation in \eqref{PDE-linearized}. We mollify only the coefficient of $\varphi_{yy}$, as the coefficients of $\varphi_x$, $\varphi_y$ pose no problems. 

\begin{lemma} \label{lem:existence-estimate}
Let $q\in L^1_{\mathrm{prob}}(\R^2)$ and $n\in\N$. Under Assumption \ref{asmp:main} there exists a unique probability measure $\mu_n$ on $\R^2$ with a finite first moment such that 
\begin{align}
\begin{split}
0 \!=\! \int_{\R^2} \! \bigg(\frac12\,\sigma_1^2(x)\,\varphi_{xx}(x,y) + \frac12\,\sigma_2^2(y)\,f^{-2}(y)\,G^{f^{-2};q}_n(x)\,\varphi_{yy}(x,y) 
\qquad\qquad \qquad\qquad\qquad\quad\;\;\; \\
+\,b_1(x)\,(hf^{-2})(y)\,G^{hf^{-2};q}(x)\,\varphi_x(x,y)
\!+\! b_2(y)\,f^{-2}(y)\,G^{f^{-2};q}(x)\,\varphi_y(x,y) 
\!\bigg)\mu_n(\mathrm{d}x,\mathrm{d}y), \\ \varphi \in C^\infty_c(\R^2).
\end{split} \label{PDE-linearized-mollified}
\end{align}
Further, there are constants $\widetilde{C}_1,\widetilde{C}_2<\infty$ depending only on the constants mentioned in Assumption \ref{asmp:main} (and, in particular, independent of $q$ and $n$), such that $\mu_n(\mathrm{d}x,\mathrm{d}y) = \widetilde{q}_n(x,y)\,\mathrm{d}x\,\mathrm{d}y$ for some $\widetilde{q}_n \in L^1_{\mathrm{prob}}(\R^2)$ satisfying $\sqrt{\widetilde{q}_n} \in W^1_2(\R^2)$ and
\begin{align}
& \frac14\int_{\R^2}\frac{|\nabla\widetilde{q}_n(x,y) |^2}{\widetilde{q}_n(x,y)}\,\mathrm{d}x\,\mathrm{d}y 
= \int_{\R^2} \big|\nabla\sqrt{\widetilde{q}}_n(x,y)\big|^2\,\mathrm{d}x\,\mathrm{d}y \le \widetilde{C}_1,  \label{le:mollified-gradbound}  \\
& \int_{\R^2} (x^2+y^2)\,\mu_n(\mathrm{d}x,\mathrm{d}y) \le \widetilde{C}_2. \label{le:mollified-momentbound}
\end{align}
\end{lemma}

\smallskip

\noindent\textbf{Proof.} By Proposition \ref{pr:general linear SDE}, there is a unique stationary weak solution of the SDE
\begin{align*}
\begin{cases}
\;\mathrm{d}\widetilde{X}^n_t = b_1(\widetilde{X}^n_t)\,(hf^{-2})(\widetilde{Y}^n_t)\,G^{hf^{-2};q}(\widetilde{X}^n_t)\,\mathrm{d}t + \sigma_1(\widetilde{X}^n_t)\,\mathrm{d}W_t, \\
\;\mathrm{d}\widetilde{Y}^n_t \,= b_2(\widetilde{Y}^n_t)\,f^{-2}(\widetilde{Y}^n_t)\,G^{f^{-2};q}(\widetilde{X}^n_t)\,\mathrm{d}t + \sigma_2(\widetilde{Y}^n_t)\,f^{-1}(\widetilde{Y}^n_t)\,\sqrt{G_n^{f^{-2};q}(\widetilde{X}^n_t)}\,\mathrm{d}B_t.
\end{cases}
\end{align*}
According to \cite[Theorem 2.5]{trevisan}, every probability measure $\mu_n$ on $\R^2$ with a finite first moment solving \eqref{PDE-linearized-mollified} can be identified with the fixed-time marginal distributions of $(\widetilde{X}^n_t,\widetilde{Y}^n_t)_{t\ge0}$. Note hereby that, for any fixed $t\ge0$, the test functions of \cite[Definition 2.2]{trevisan} belong to $C_c^2(\R^2)$ and the equation in \eqref{PDE-linearized-mollified} holds for such functions due to a straightforward density argument. In particular, thanks to Proposition \ref{pr:general linear SDE}, 
\begin{align}\label{2nd_moment}
\int_{\R^2} (x^2+y^2)\,\mu_n(\mathrm{d}x,\mathrm{d}y) 
= \E\big[(\widetilde{X}^n_t)^2+(\widetilde{Y}^n_t)^2\big]
\le \widetilde{C}_2,
\end{align}
 where the constant $\widetilde{C}_2 <\infty$ depends only on the constants mentioned in Assumption \ref{asmp:main} (and not on $q$ or $n$).

\medskip

It remains to show the existence of a density $\widetilde{q}_n$ with the properties described in the lemma. Define $A(x,y) = (A^{ij}(x,y))_{1\le i,j\le 2}$ and $B(x,y)=(B^i(x,y))_{1\le i\le 2}$ by
\begin{align*}
A(x,y) \!=\! \begin{pmatrix}
\frac12\,\sigma_1^2(x) & 0 \\
0 & \frac12\,\sigma_2^2(y)\,f^{-2}(y)\,G_n^{f^{-2};q}(x)
\end{pmatrix},\;\;\;
B(x,y) \!=\! \begin{pmatrix}
b_1(x)\,(hf^{-2})(y)\,G^{hf^{-2};q}(x) \\
b_2(y)\,f^{-2}(y)\,G^{f^{-2};q}(x)
\end{pmatrix}.
\end{align*}
Define also $D(x,y)=(D^i(x,y))_{1\le i\le 2}$ by
\begin{align*}
D^1(x,y) &= A^{11}_x(x,y) + A^{21}_y(x,y) = \sigma_1(x)\,\sigma_1'(x), \\
D^2(x,y) &= A^{12}_x(x,y) + A^{22}_y(x,y) = \big(\sigma_2(y)\,\sigma_2'(y)\,f^{-2}(y)
-\sigma_2^2(y)\,f^{-3}(y)\,f'(y)\big)\,G_n^{f^{-2};q}(x).
\end{align*}
In view of Assumption \ref{asmp:main} and the inequality in \eqref{def:mollified-gradbound}, the function $A$ is uniformly Lipschitz continuous. Moreover, $A$ is bounded and its eigenvalues are bounded below by $\alpha=\min(\inf_\R \sigma_1^2,\inf_\R (\sigma_2^2 f^{-2})\,\inf_\R f^2)/2>0$. In addition, with the notation $\|\cdot\|_\infty$ for the supremum norm, Assumption \ref{asmp:main} yields
\begin{align}
\begin{split}
|B^1(x,y)| &\le C_2(1+|x|)\,\|hf^{-2}\|_\infty\,\|h^{-1}f^2\|_\infty, \\
|B^2(x,y)| &\le C_2(1+|y|)\,\|f^{-2}\|_\infty\,\|f^2\|_\infty.
\end{split} \label{pf:driftbound1}
\end{align}
Thus, the estimate \eqref{2nd_moment} renders Proposition \ref{BKR_reg} applicable. We conclude that $\mu_n$ is of the form $\mu_n(dx,dy)=\widetilde{q}_n(x,y)\,\mathrm{d}x\,\mathrm{d}y$, for some $\widetilde{q}_n \in L^1_{\mathrm{prob}}(\R^2)$ with $\sqrt{\widetilde{q}_n} \in W^1_2(\R^2)$ and
\begin{align}\label{BKRbound}
\frac14\int_{\R^2}\frac{|\nabla\widetilde{q}_n(x,y) |^2}{\widetilde{q}_n(x,y)}\,\mathrm{d}x\,\mathrm{d}y 
= \int_{\R^2} \big|\nabla\sqrt{\widetilde{q}}_n(x,y)\big|^2\,\mathrm{d}x\,\mathrm{d}y \le \frac{1}{2\alpha^2}\int_{\R^2} \big(|B|^2+|D|^2\big)\,\mathrm{d}\mu_n.
\end{align}
By Assumption \ref{asmp:main}, we have
\begin{align*}
& \|D^1\|_\infty \le \|\sigma_1\|_\infty\,\|\sigma'_1\|_\infty<\infty, \\
& \|D^2\|_\infty \le \big(\|\sigma_2\|_\infty\,\|\sigma'_2\|_\infty\,\|f^{-2}\|_\infty
+\|\sigma_2^2\|_\infty\,\|f^{-3}\|_\infty\,\|f'\|_\infty\big)\,\|f^2\|_\infty<\infty,
\end{align*}
which together with \eqref{pf:driftbound1} and \eqref{2nd_moment} allows us to bound the rightmost expression in \eqref{BKRbound} by a constant $\widetilde{C}_1$ as in the statement of the lemma. \ep

\medskip

With Lemma \ref{lem:existence-estimate} now established, we are going to send $n\to\infty$ to prove Proposition \ref{pr:linearized-PDE-estimates}. In order to facilitate this, we prepare a compactness lemma that is used again below. With the constants $\widetilde{C}_1,\widetilde{C}_2$ of Lemma \ref{lem:existence-estimate}, let
\begin{align*}
K_0=\bigg\{\!q\in L^1_{\mathrm{prob}}(\R^2):\, \sqrt{q} \in W^1_2(\R^2),\,\int_{\R^2} \! |\nabla \sqrt{q}|^2\,\mathrm{d}x\,\mathrm{d}y \le \widetilde{C}_1,\,\int_{\R^2} \! (x^2+y^2)\,q\,\mathrm{d}x\,\mathrm{d}y \le \widetilde{C}_2\bigg\}.
\end{align*}

\begin{lemma} \label{lem:K0-compact}
$K_0$ is norm-compact in $L^1(\R^2)$.
\end{lemma}

\smallskip

\noindent\textbf{Proof.} We first show that $K_0$ is norm-precompact. For any sequence $(q_n)_{n\in\N}$ in $K_0$, by the Cauchy-Schwarz inequality, 
\begin{align*}
\int_{\R^2} |\nabla q_n|\,\mathrm{d}x\,\mathrm{d}y
& = \int_{\R^2} 2\sqrt{q_n}\,|\nabla\sqrt{q_n}|\,\mathrm{d}x\,\mathrm{d}y \\
& \le 2\bigg(\int_{\R^2} q_n\,\mathrm{d}x\,\mathrm{d}y\bigg)^{1/2} 
\bigg(\int_{\R^2} |\nabla\sqrt{q_n}|^2\,\mathrm{d}x\,\mathrm{d}y \bigg)^{1/2} 
\le 2\sqrt{\widetilde{C}_1},\quad n\in\N.
\end{align*}
The Rellich-Kondrachov theorem (see, e.g., \cite[Theorem 5.7.1]{evans10}), employed on the open balls $\{(x,y)\in\R^2:\,x^2+y^2<N^2\}$, $N\in\N$ and followed by a diagonalization argument, gives rise to a subsequence of $(q_n)_{n\in\N}$, also referred to as $(q_n)_{n\in\N}$, converging locally in $L^1(\R^2)$ and a.e.\ to some locally integrable $q$. Thanks to Markov's inequality,
\begin{align*}
\int_{\{(x,y)\in\R^2:\,x^2+y^2\ge N^2\}} q_n\,\mathrm{d}x\,\mathrm{d}y 
\le \widetilde{C}_2 N^{-2}, \quad N,n\in\N.
\end{align*}
This and Fatou's lemma imply
\begin{align*}
\int_{\R^2}|q_n - q|\,\mathrm{d}x\,\mathrm{d}y 
\le \int_{\{(x,y)\in\R^2:\,x^2+y^2< N^2\}} |q_n - q|\,\mathrm{d}x\,\mathrm{d}y  
+ 2\widetilde{C}_2N^{-2}, \quad N,n\in\N.
\end{align*}
Taking $n\to\infty$ and then $N\to\infty$ we deduce that $q_n\to q$ in $L^1(\R^2)$ as $n\to\infty$. 

\medskip

It remains to check that $K_0$ is closed. Let $(q_n)_{n\in\N}$ be a sequence in $K_0$ that tends in $L^1(\R^2)$ to some $q$. By passing to an a.e.\ convergent subsequence and applying Fatou's lemma we get $\int_{\R^2} (x^2+y^2)\,q(x,y)\,\mathrm{d}x\,\mathrm{d}y \le \widetilde{C}_2$. Finally, from the bound
\begin{align*}
\int_{\R^2} \sqrt{q_n}^2\,\mathrm{d}x\,\mathrm{d}y
+\int_{\R^2} |\nabla\sqrt{q_n}|^2\,\mathrm{d}x\,\mathrm{d}y 
\le 1+\widetilde{C}_1,\quad n\in\N
\end{align*}
we conclude that there exists a subsequence of $(\sqrt{q_n})_{n\in\N}$ tending weakly in $W^1_2(\R^2)$ and locally strongly in $L^2(\R^2)$ to some $u$ (cf.\ \cite[Chapter 10, Theorem 7]{lax} and \cite[remark on p.~274]{evans10}). Moreover, because $q_n \to q$ in $L^1(\R^2)$, we must have that $q=u^2$ (the locally strong $L^1(\R^2)$-limit of the subsequence). Thus, $(\sqrt{q_n})_{n\in\N}$ converges weakly in $W^1_2(\R^2)$ to $\sqrt{q}$. By the weak lower semi-continuity of the $W^1_2(\R^2)$-norm (see, e.g., \cite[Chapter 10, Theorem 5]{lax}), 
\begin{align*}
1+\int_{\R^2} |\nabla \sqrt{q}|^2\,\mathrm{d}x\,\mathrm{d}y 
& = \int_{\R^2} \sqrt{q}^2\,\mathrm{d}x\,\mathrm{d}y 
+\int_{\R^2} |\nabla \sqrt{q}|^2\,\mathrm{d}x\,\mathrm{d}y \\
& \le \liminf_{n\to\infty} \bigg(\int_{\R^2} \sqrt{q_n}^2\,\mathrm{d}x\,\mathrm{d}y 
+\int_{\R^2} |\nabla \sqrt{q_n}|^2\,\mathrm{d}x\,\mathrm{d}y\bigg) 
\le 1+\widetilde{C}_1,
\end{align*}
which finishes the proof. \ep

\medskip

\noindent\textbf{Proof of Proposition \ref{pr:linearized-PDE-estimates}.} The existence and uniqueness of a probability measure $\mu$ on $\R^2$ with a finite first moment satisfying \eqref{PDE-linearized} and the claims (ii) and (iii) follow from Proposition \ref{pr:general linear SDE} and \cite[Theorem 2.5]{trevisan}, as in the beginning of the proof of Lemma \ref{lem:existence-estimate}. To obtain claim (i) we recall the densities $(\widetilde{q}_n)_{n \in \N}$ of Lemma \ref{lem:existence-estimate}, and we aim to find a subsequential limit $\widetilde{q}$ of this sequence in $L^1(\R^2)$ which enjoys the properties in (i) and such that $\widetilde{q}(x,y)\,\mathrm{d}x\,\mathrm{d}y$ has a finite first moment and solves \eqref{PDE-linearized}. Then, by uniqueness,  $\mu(\mathrm{d}x,\mathrm{d}y)=\widetilde{q}(x,y)\,\mathrm{d}x\,\mathrm{d}y$.

\medskip

Since $\widetilde{q}_n \in K_0$ for all $n$, Lemma \ref{lem:K0-compact} allows us to extract a subsequence converging in $L^1(\R^2)$ to some $\widetilde{q}\in K_0$. To verify that $\widetilde{q}(x,y)\,\mathrm{d}x\,\mathrm{d}y$ solves \eqref{PDE-linearized} we set, for $\varphi \in C^\infty_c(\R^2)$,
\begin{align*}
({\mathcal L}^q_n\varphi)(x,y)=&\;\frac12\,\sigma_1^2(x)\,\varphi_{xx}(x,y) 
+ \frac12\,\sigma_2^2(y)\,f^{-2}(y)\,G^{f^{-2};q}_n(x)\,\varphi_{yy}(x,y) \\
& + b_1(x)\,(hf^{-2})(y)\,G^{hf^{-2};q}(x)\,\varphi_x(x,y)
+b_2(y)\,f^{-2}(y)\,G^{f^{-2};q}(x)\,\varphi_y(x,y).
\end{align*}
Note that $G_n^{f^{-2};q} \to G^{f^{-2};q}$ a.e.\ when $n\to\infty$ by the Lebesgue differentiation theorem, and the functions ${\mathcal L}^q_n\varphi$ are bounded uniformly in $n\in\N$. Thus, with  
\begin{align*}
({\mathcal L}^q\varphi)(x,y)= &\;\frac12\,\sigma_1^2(x)\,\varphi_{xx}(x,y) 
+ \frac12\,\sigma_2^2(y)\,f^{-2}(y)\,G^{f^{-2};q}(x)\,\varphi_{yy}(x,y) \\
&+ b_1(x)\,(hf^{-2})(y)\,G^{hf^{-2};q}(x)\,\varphi_x(x,y)
+ b_2(y)\,f^{-2}(y)\,G^{f^{-2};q}(x)\,\varphi_y(x,y),
\end{align*}
we find that
\begin{align*}
\bigg|\int_{\R^2} ({\mathcal L^q}\varphi)\,\widetilde{q}\,\mathrm{d}x\,\mathrm{d}y \!-\!\int_{\R^2} ({\mathcal L}^q_n\varphi)\,\widetilde{q}_n\,\mathrm{d}x\,\mathrm{d}y \bigg| 
\le \int_{\R^2} |{\mathcal L^q}\varphi\!-\!{\mathcal L}^q_n\varphi|\,\widetilde{q}\,\mathrm{d}x\,\mathrm{d}y
\!+\! \int_{\R^2} |{\mathcal L}^q_n\varphi|\,|\widetilde{q}\!-\!\widetilde{q}_n|\,\mathrm{d}x\,\mathrm{d}y
\end{align*}
tends to $0$ as $n\to\infty$. Consequently, 
\begin{align*}
\int_{\R^2} ({\mathcal L^q}\varphi)\,\widetilde{q}\,\mathrm{d}x\,\mathrm{d}y = \lim_{n\to\infty} \int_{\R^2} ({\mathcal L}^q_n\varphi)\,\widetilde{q}_n\,\mathrm{d}x\,\mathrm{d}y = 0,
\end{align*}
i.e., the probability measure $\widetilde{q}\,\mathrm{d}x\,\mathrm{d}y$ (which has a finite first moment) solves \eqref{PDE-linearized}. \ep


\section{Proof of Theorem \ref{th:intro:main}} \label{se:2.2}

\subsection{Continuity of the conditional expectation}

Our proof of Theorem \ref{th:intro:main} proceeds via a fixed-point argument. The continuity of the underlying fixed-point map relies on the next lemma, which we prepare beforehand. 

\begin{lemma} \label{lem:G[]-continuous}
Let $\psi:\,\R\to(0,\infty)$ be a measurable function with $c_\psi\le\psi\le C_\psi$ for some $c_\psi,C_\psi\in(0,\infty)$. Suppose $(q_n)_{n\in\N}$ is a sequence in $L^1_{\mathrm{prob}}(\R^2)$ converging to some $q\in L^1_{\mathrm{prob}}(\R^2)$. Then, with $m(x)=\int_\R q(x,y)\,\mathrm{d}y$, we have $G^{\psi;q_n}\to G^{\psi;q}$ in $L^1(\R,m(x)\,\mathrm{d}x)$.
\end{lemma}

\smallskip

\noindent\textbf{Proof.} Define the marginal densities $m_n(x)= \int_{\R} q_n(x,y)\,\mathrm{d}y$, $n\in\N$. By using the uniform Lipschitz continuity of $x\mapsto 1/x$ on $(c_\psi,\infty)$ and applying the triangle inequality repeatedly we derive the estimates
\begin{align*}
&\,\int_\R |G^{\psi;q_n}(x) - G^{\psi;q}(x)|\,m(x)\,\mathrm{d}x \\
&= \int_\R \bigg|\frac{m_n(x)}{\int_\R \psi(y)\,q_n(x,y)\,\mathrm{d}y} 
- \frac{m(x)}{\int_\R \psi(y)\,q(x,y)\,\mathrm{d}y}\bigg|\,m(x)\,\mathrm{d}x \\
&\le (c_\psi)^{-2} \int_\R \bigg|\frac{\int_\R \psi(y)\,q_n(x,y)\,\mathrm{d}y}{m_n(x)} - \frac{\int_\R \psi(y)\,q(x,y)\,\mathrm{d}y}{m(x)}\bigg|\,m(x)\,\mathrm{d}x \\
&\le (c_\psi)^{-2}\int_\R \frac{\int_\R \psi(y)\,q_n(x,y)\,\mathrm{d}y}{m_n(x)}\,|m(x)-m_n(x)|\,\mathrm{d}x \\
&\quad\, + (c_\psi)^{-2} \int_\R \bigg|\int_\R \psi(y)\,q_n(x,y)\,\mathrm{d}y - \int_\R \psi(y)\,q(x,y)\,\mathrm{d}y\bigg|\,\mathrm{d}x \\
&\le (c_\psi)^{-2}C_\psi\int_\R \frac{\int_\R q_n(x,y)\,\mathrm{d}y}{m_n(x)}\,|m(x)-m_n(x)|\,\mathrm{d}x 
+ (c_\psi)^{-2}C_\psi\int_{\R^2} |q_n(x,y)-q(x,y)|\,\mathrm{d}x\,\mathrm{d}y \\
&= (c_\psi)^{-2}C_\psi \int_\R |m(x)-m_n(x)|\,\mathrm{d}x 
+(c_\psi)^{-2}C_\psi \int_{\R^2} |q_n(x,y)-q(x,y)|\,\mathrm{d}x\,\mathrm{d}y \\
&\le 2(c_\psi)^{-2}C_\psi\int_{\R^2} |q_n(x,y)-q(x,y)|\,\mathrm{d}x\,\mathrm{d}y.
\end{align*}
The latter expression tends to $0$ as $n\to\infty$ by assumption. \ep

\subsection{Main line of the argument}

We are now ready to present the main line of the argument. With the constants $\widetilde{C}_1$, $\widetilde{C}_2$ and $\widetilde{\delta}_R$, $R\in(0,\infty)$ of Proposition \ref{pr:linearized-PDE-estimates}, we set
\begin{align*}
K=\bigg\{\!q\in L^1_{\mathrm{prob}}(\R^2)\!:\sqrt{q}\in W^1_2(\R^2),\,\int_{\R^2} |\nabla \sqrt{q}|^2\,\mathrm{d}x\,\mathrm{d}y\le\widetilde{C}_1,\,\int_{\R^2} (x^2+y^2)\,q\,\mathrm{d}x\,\mathrm{d}y\le\widetilde{C}_2,\quad \\
\forall R\in(0,\infty):\,\int_\R q(x,y)\,\mathrm{d}y\ge\widetilde{\delta}_R\text{ a.e. in }[-R,R]\bigg\}.
\end{align*}
In other words, $K$ is the intersection of $K_0$ (introduced prior to Lemma \ref{lem:K0-compact}) with the set 
\begin{align*}
K_1=\bigg\{q\in L^1_{\mathrm{prob}}(\R^2):\,\int_\R q(x,y)\,\mathrm{d}y\ge\widetilde{\delta}_R\text{ a.e. in }[-R,R],\;R\in(0,\infty)\bigg\}.
\end{align*}
In addition, we let $\Phi:\,K\to K$ be the map taking each $q \in K$ to the unique density $\widetilde{q}\in K$ given in Proposition \ref{pr:linearized-PDE-estimates}. The following two lemmas establish some basic properties of $K$ and $\Phi$.  

\begin{lemma} \label{lem:Phi-compact}
$K$ is norm-compact and convex in $L^1(\R^2)$. 
\end{lemma}

\begin{lemma} \label{lem:Phi-continuous}
The map $\Phi:\,K\to K$ is continuous with respect to the $L^1(\R^2)$--norm.
\end{lemma}

\noindent\textbf{Proof of Lemma \ref{lem:Phi-compact}.} As noted above, $K=K_0 \cap K_1$. Since $K_0$ is norm-compact by Lemma \ref{lem:K0-compact}, it is enough to show that $K_1$ is closed. But this is straightforward: If a sequence $(q_n)_{n\in\N}$ in $K_1$ converges to some $q\in L^1_{\mathrm{prob}}(\R^2)$, then we have, for every $R>0$ and every measurable $\psi:\,[-R,R] \to [0,1]$, 
\begin{align*}
\int_{-R}^R \psi(x)\,\int_{\R} q_n(x,y)\,\mathrm{d}y\,\mathrm{d}x \ge 
\widetilde{\delta}_R\,\int_{-R}^R\psi(x)\,\mathrm{d}x,\quad n\in\N. 
\end{align*}
Passing to the limit $n\to\infty$ we get the same inequality with $q_n$ replaced by $q$, which implies $\int_\R q(x,y)\,\mathrm{d}y \ge \widetilde{\delta}_R$ a.e.\ in $[-R,R]$, for every $R > 0$. 

\medskip

It remains to prove the convexity of $K$. Clearly, the only delicate point is the stability of the properties $\sqrt{q}\in W^1_2(\R^2)$ and $\int_{\R^2} |\nabla \sqrt{q}|^2\,\mathrm{d}x\,\mathrm{d}y\le\widetilde{C}_1$ under convex combinations. This fact is fairly well-known, as the functional $\frac12\int_{\R^2} |\nabla \sqrt{q}|^2\,\mathrm{d}x\,\mathrm{d}y \!=\! \frac18\int_{\R^2} |\nabla\log q|^2\, q\,\mathrm{d}x\,\mathrm{d}y$, often called the Fisher information, provides the rate function in the work of Donsker and Varadhan \cite{donsker1975asymptotic} on the large deviations for the occupation measure of Brownian motion. Nonetheless, we describe a short self-contained proof: Observe that $\sqrt{q}\in W^1_2(\R^2)$ if and only if $q\in W^1_1(\R^2)$ and $|\nabla q|^2/q \in L^1(\R^2)$. For such $q$ we may write
\begin{align*}
\int_{\R^2} |\nabla \sqrt{q}|^2 \,\mathrm{d}x\,\mathrm{d}y &= \frac14 \int_{\R^2} |\nabla \log q|^2\,q\,\mathrm{d}x\,\mathrm{d}y \\
&= \int_{\R^2} \sup_{z\in\R^2} \big(z\cdot\nabla\log q(x,y)-|z|^2\big)\,q(x,y)\,\mathrm{d}x\,\mathrm{d}y \\
&= \sup_{\eta \in L^\infty(\R^2)}\,\int_{\R^2} \big(\eta(x,y)\cdot\nabla\log q(x,y)-|\eta(x,y)|^2\big)\,q(x,y)\,\mathrm{d}x\,\mathrm{d}y \\
&= \sup_{\eta \in L^\infty(\R^2)}\,\int_{\R^2} \big(\eta(x,y)\cdot\nabla q(x,y) - |\eta(x,y)|^2\,q(x,y)\big)\,\mathrm{d}x\,\mathrm{d}y.
\end{align*}
As a supremum of linear functionals, $\!\int_{\R^2}\!|\nabla\!\sqrt{q}|^2 \,\mathrm{d}x\,\mathrm{d}y$ is convex on $L^1_{\mathrm{prob}}(\R^2)\cap W^1_1(\R^2)$. \ep

\medskip

\noindent\textbf{Proof of Lemma \ref{lem:Phi-continuous}.} Suppose $(q_n)_{n\in\N}$ is a convergent sequence in $K$ with a limit $q\in K$. Define $\widetilde{q}_n=\Phi(q_n)\in K$ for $n\in\N$. By Lemma \ref{lem:Phi-compact}, any subsequence of $(\widetilde{q}_n)_{n\in\N}$ has an $L^1(\R^2)$--convergent subsubsequence, and we aim to verify that any resulting limit point $\widetilde{q}\in K$ equals  $\Phi(q)$. To this end, we relabel the subsubsequence so that $\widetilde{q}_n\to\widetilde{q}$ in $L^1(\R^2)$. The definition of $\Phi$ yields
\begin{align}
\begin{split}
\int_{\R^2} ({\mathcal L}^{q_n}\varphi)\,\widetilde{q}_n\,\mathrm{d}x\,\mathrm{d}y
=0, \quad \varphi \in C^\infty_c(\R^2),
\end{split}  \label{pf:PDEatlimit1}
\end{align}
where we have employed the notation
\begin{align*}
({\mathcal L}^{q_n}\varphi)(x,y) = &\;\frac12\,\sigma_1^2(x)\,\varphi_{xx}(x,y) + \frac12\,\sigma_2^2(y)\,f^{-2}(y)\,G^{f^{-2};q_n}(x)\,\varphi_{yy}(x,y) \\
& + b_1(x)\,(hf^{-2})(y)\,G^{hf^{-2};q_n}(x)\,\varphi_x(x,y) 
+ b_2(y)\,f^{-2}(y)\,G^{f^{-2};q_n}(x)\,\varphi_y(x,y).
\end{align*}
With $m(x) = \int_\R q(x,y)\,\mathrm{d}y$, Lemma \ref{lem:G[]-continuous} gives $G^{f^{-2};q_n}\to G^{f^{-2};q}$ and $G^{hf^{-2};q_n}\to G^{hf^{-2};q}$ in $L^1(\R,m(x)\,\mathrm{d}x)$. By taking a further subsequence, we can ensure that both convergences hold also a.e.\ with respect to the probability measure $m(x)\,\mathrm{d}x$. Since $q\in K$, we have $m>0$ (Lebesgue) a.e., and so $G^{f^{-2};q_n}\to G^{f^{-2};q}$ and $G^{hf^{-2};q_n}\to G^{hf^{-2};q}$ (Lebesgue) a.e.\ along the same subsequence. Thus, for all $\varphi \in C^\infty_c(\R^2)$,
\begin{align*}
\bigg|\!\int_{\R^2} ({\mathcal L^q}\varphi)\,\widetilde{q}\,\mathrm{d}x\,\mathrm{d}y \!-\!\int_{\R^2} ({\mathcal L}^{q_n}\varphi)\,\widetilde{q}_n\,\mathrm{d}x\,\mathrm{d}y \bigg| 
\!\le\! \int_{\R^2}\! |{\mathcal L^q}\varphi\!-\!{\mathcal L}^{q_n}\varphi|\,\widetilde{q}\,\mathrm{d}x\,\mathrm{d}y
\!+\! \int_{\R^2}\! |{\mathcal L}^{q_n}\varphi|\,|\widetilde{q}\!-\!\widetilde{q}_n|\,\mathrm{d}x\,\mathrm{d}y
\end{align*}
tends to $0$ along that subsequence. The definition of $\Phi(q)$ now shows that $\widetilde{q} = \Phi(q)$. \ep

\medskip

\noindent\textbf{Proof of Theorem \ref{th:intro:main}.} Lemmas \ref{lem:Phi-compact} and \ref{lem:Phi-continuous} reveal that $\Phi$ is a continuous map from the (non-empty) norm-compact and convex set $K\subset L^1(\R^2)$ into itself. Consequently, there exists a fixed-point $\widetilde{p}=\Phi(\widetilde{p})\in K$ by the Schauder fixed-point theorem. In view of the definition of $\Phi$, the latter satisfies the transformed stationary Fokker-Planck equation \eqref{mainpde_weak-transformed} for all $\varphi\in C^\infty_c(\R^2)$. Parts (ii) and (iv) of Proposition \ref{prop_changemeas} then imply that $p(x,y)=\widetilde{p}(x,y)\,f^{-2}(y)\,G^{f^{-2};\widetilde{p}}(x)$ satisfies the original stationary Fokker-Planck equation \eqref{mainpde_weak} for all $\varphi\in C^\infty_c(\R^2)$. Moreover, there is a unique stationary weak solution of the SDE
\begin{align*}
\begin{cases}
\;\mathrm{d}X_t = b_1(X_t)\,h(Y_t)\,G^{h;p}(X_t)\,\mathrm{d}t + \sigma_1(X_t)\,f(Y_t)\,\sqrt{G^{f^2;p}(X_t)}\,\mathrm{d}W_t, \\
\;\mathrm{d}Y_t \,=\, b_2(Y_t)\,\mathrm{d}t + \sigma_2(Y_t)\,\mathrm{d}B_t
\end{cases} 
\end{align*}
by Proposition \ref{pr:general linear SDE}. Applying \cite[Theorem 2.5]{trevisan} as in the beginning of the proof of Lemma \ref{lem:existence-estimate} (noting that $p(x,y)\,\mathrm{d}x\,\mathrm{d}y$ has a finite first moment because $\widetilde{p}(x,y)\,\mathrm{d}x\,\mathrm{d}y$ does), the distributions  ${\mathcal L}(X_t,Y_t)$ for $t\ge0$ are readily identified with $p(x,y)\,\mathrm{d}x\,\mathrm{d}y$. Hence, $(X_t,Y_t)_{t\ge0}$ is a weak solution of the SDE \eqref{intro:mainSDEsystem} which obeys $(X_t,Y_t)\stackrel{d}{=}(X_0,Y_0)$ for $t\ge0$.  

\medskip

Claims (i) and (ii) are immediate corollaries of Proposition \ref{th:intro:Markovprojection} and Lemma \ref{le:intro:onedimensional}, and we turn to the proof of claim (iii). Inserting the definition of $G^{f^{-2};\widetilde{p}}$ into $p(x,y)=\widetilde{p}(x,y)\,f^{-2}(y)\,G^{f^{-2};\widetilde{p}}(x)$, differentiating via the product rule, and using the boundedness of $f^{-2}$, $G^{f^{-2};\widetilde{p}}$, $\widetilde{p}/p$ and $f^2$ we estimate $\int_{\R^2} \frac{p_x^2}{p}\,\mathrm{d}x\,\mathrm{d}y$ by a constant multiple of 
\begin{align*}
\int_{\R^2} \frac{\widetilde{p}_x^2}{\widetilde{p}}\,\mathrm{d}x\,\mathrm{d}y
+\int_{\R^2} \widetilde{p}\;\frac{(\int_\R \widetilde{p}_x(\cdot,z)\,\mathrm{d}z)^2}{(\int_\R \widetilde{p}(\cdot,z)\,\mathrm{d}z)^2}\,\mathrm{d}x\,\mathrm{d}y
+\int_{\R^2} \widetilde{p}\;\frac{(\int_\R f^{-2}(z)\,\widetilde{p}_x(\cdot,z)\,\mathrm{d}z)^2}{(\int_\R  \widetilde{p}(\cdot,z)\,\mathrm{d}z)^2}\,\mathrm{d}x\,\mathrm{d}y.
\end{align*}
The first $\mathrm{d}x\,\mathrm{d}y$--integral is finite thanks to $\widetilde{p}\in K$. Setting $\psi=1$ and $\psi=f^{-2}$ in the cases of the second and third $\mathrm{d}x\,\mathrm{d}y$--integrals, respectively, and integrating in $y$ we end up with
\begin{align*}
\int_\R \bigg(\int_\R
\psi(z)\,\widetilde{p}_x(x,z)\,\mathrm{d}z\bigg)^2\,\frac{1}{\widetilde{m}_1(x)}\,\mathrm{d}x,
\end{align*} 
where $\widetilde{m}_1(x)=\int_\R  \widetilde{p}(x,z)\,\mathrm{d}z$. By Jensen's inequality, 
\begin{equation}\label{Dan's Jensen}
\begin{split}
\int_\R \bigg(\int_\R \psi(z)\,\widetilde{p}_x(x,z)\,\mathrm{d}z\bigg)^2\,\frac{1}{\widetilde{m}_1(x)}\,\mathrm{d}x
& =\int_\R \bigg(\int_\R \psi(z)\,\frac{\widetilde{p}_x(x,z)}{\widetilde{p}(x,z)}\,\frac{\widetilde{p}(x,z)}{\widetilde{m}_1(x)}\,\mathrm{d}z\bigg)^2\,\widetilde{m}_1(x)\,\mathrm{d}x \\
& \le \int_\R \int_\R \psi^2(z)\,\frac{\widetilde{p}_x(x,z)^2}{\widetilde{p}(x,z)^2}\,\frac{\widetilde{p}(x,z)}{\widetilde{m}_1(x)}\,\mathrm{d}z\,\widetilde{m}_1(x)\,\mathrm{d}x,
\end{split}
\end{equation}
which is finite due to the boundedness of $\psi$ and $\widetilde{p}\in K$. It follows that $\int_{\R^2} \frac{p_x^2}{p}\,\mathrm{d}x\,\mathrm{d}y<\infty$. Differentiating by means of the product rule and relying on the boundedness of $\widetilde{p}/p$, $f^{-2}$, $G^{f^{-2};\widetilde{p}}$ and $(f^{-2})'$ we see that $\int_{\R^2} \frac{p_y^2}{p}\,\mathrm{d}x\,\mathrm{d}y$ cannot exceed a constant multiple of $\int_{\R^2} \frac{\widetilde{p}_y^2}{\widetilde{p}}+\widetilde{p}\,\mathrm{d}x\,\mathrm{d}y$. In view of $\widetilde{p}\in K$, we obtain claim (iii).  

\medskip

Lastly, we prove that any stationary weak solution of the SDE \eqref{intro:mainSDEsystem} must be strong. For this purpose, we recall a theorem of Veretennikov (cf.\ \cite[Theorem 4]{veretennikov1980strong}): Strong existence and pathwise uniqueness are valid for a one-dimensional SDE 
\begin{align*}
\mathrm{d}Z_t = b(t,Z_t)\,\mathrm{d}t + \sigma(t,Z_t)\,\mathrm{d}\beta_t 
\end{align*}
under the assumptions that $b$ and $\sigma$ are bounded and measurable, $\inf_{t\ge0,\,z\in\R} \sigma(t,z) > 0$ and $\sup_{t\ge0} |\sigma(t,z)-\sigma(t,\widetilde{z})|\le C|z-\widetilde{z}|^{1/2}$ for all $z,\widetilde{z}\in\R$, with some $C<\infty$. By a straightforward localization argument, one can extend his result to drifts fulfilling the linear growth condition $|b(t,z)|\le C(1+|z|)$, $t\ge0$, $z\in\R$, as well as relax the regularity condition on the diffusion coefficient to local $\tfrac12$--H\"older continuity, in the sense that for each $R >0$ there is a $C_R<\infty$ with $\sup_{t\ge 0}|\sigma(t,z)-\sigma(t,\widetilde{z})| \le C_R|z-\widetilde{z}|^{1/2}$ for all $z,\widetilde{z}\in[-R,R]$. 

\medskip

Next, suppose $(X_t,Y_t)_{t\ge0}$ is a stationary weak solution of the SDE \eqref{intro:mainSDEsystem}, defined on some probability space $(\Omega,\F,\FF,\PP)$ supporting independent standard $\FF$-Brownian motions $W$ and $B$, so that $(X_t)_{t\ge0}$ and $(Y_t)_{t\ge0}$ are of course adapted to the same filtration $\FF$. Strong existence and pathwise uniqueness hold for the one-dimensional SDE $(Y_t)_{t\ge0}$ solves, and so each $Y_t$ is measurable with respect to $(Y_0,(B_s)_{s\in[0,t]})$. In addition, the processes $(W_t)_{t\ge0}$ and $(Y_t)_{t\ge0}$ are conditionally independent given $X_0$, because the Brownian motion $(W_t)_{t\ge0}$ is independent of $(X_0,Y_0)$ and the Brownian motion $(B_t)_{t\ge0}$. Thus, for any $x_0\in\R$ and any continuous real-valued function $\mathbf{y}=(y_t)_{t\ge0}$, the conditional law of the process $(X_t)_{t\ge0}$ given $\{X_0=x_0,(Y_t)_{t\ge0}=(y_t)_{t\ge0}\}$ agrees with the law of a weak solution to the SDE
\begin{align}
\mathrm{d}X^{x_0,\mathbf{y}}_t = \lambda(X^{x_0,\mathbf{y}}_t,y_t)\,\mathrm{d}t + \xi(X^{x_0,\mathbf{y}}_t,y_t)\,\mathrm{d}\beta_t, \quad X^{x_0,\mathbf{y}}_0=x_0, \label{SDE:pathwise1}
\end{align}
living perhaps on a different probability space, where
\begin{align*}
\lambda(x,y) = b_1(x)\,\frac{h(y)}{\E[h(Y_t)|X_t=x]}, \qquad \xi(x,y)= \sigma_1(x)\,\frac{f(y)}{\sqrt{\E[f^2(Y_t)|X_t=x]}}.
\end{align*} 
We are going to show that the solution of the SDE \eqref{SDE:pathwise1} is pathwise unique. This, in turn, implies that each $X_t^{x_0,\mathbf{y}}$ is measurable with respect to $(\beta_s)_{s\in[0,t]}$. Since this is true for all $x_0$ and $\mathbf{y}$, we deduce that each $X_t$ must be measurable with respect to $(X_0,(Y_s)_{s \ge 0},(W_s)_{s\in[0,t]})$ on the original probability space. As $(Y_s)_{s\ge0}$ is adapted to $(Y_0+B_s)_{s\ge0}$, we conclude that each $X_t$ is measurable with respect to $(X_0,Y_0,(B_s)_{s \ge 0},(W_s)_{s\in[0,t]})$. Finally, the Brownian motion $(B_s-B_t)_{s \ge t}$ is independent of $\F_t$, whereas $X_t$ is $\F_t$-measurable, so that $X_t$ is, in fact, measurable with respect to $(X_0,Y_0,(B_s)_{s\in[0,t]},(W_s)_{s\in[0,t]})$.

\medskip

It remains to establish the pathwise uniqueness for the SDE \eqref{SDE:pathwise1}. Thanks to the discussion of Veretennikov's theorem above, the assumption $|b_1(x)|\le C_2(1+|x|)$ for $x\in\R$, and the boundedness of $h$ from above and below by positive constants, it suffices to check that $(t,x)\mapsto\xi(x,y_t)$ is locally $\frac{1}{2}$--H\"older continuous in $x$ uniformly in $t\ge0$, for all continuous real-valued $(y_t)_{t\ge0}$. In view of $\xi(x,y_t)=\sigma_1(x)\,f(y_t)\,\sqrt{G^{f^2;p}(x)}$, the uniform Lipschitz continuity and boundedness of $\sigma_1$, the boundedness of $f$, and the boundedness of $G^{f^2;p}$ above and below by positive constants, it is enough to prove the local $\frac{1}{2}$--H\"older continuity of $G^{f^2;p}$. To this end, we compute
\begin{align*}
(G^{f^2;p})'(x)=\frac{\int_\R p_x(x,z)\,\mathrm{d}z}{\int_\R f^2(z)\,p(x,z)\,\mathrm{d}z}-\frac{\int_\R p(x,z)\,\mathrm{d}z\,\int_\R f^2(z)\,p_x(x,z)\,\mathrm{d}z}{(\int_\R f^2(z)\,p(x,z)\,\mathrm{d}z)^2}.
\end{align*}
Hence, with $m_1(x)=\int_\R p(x,z)\,\mathrm{d}z$, we can use the lower boundedness of $f^2$ by a positive constant to estimate $\int_\R (G^{f^2;p})'(x)^2\,m_1(x)\,\mathrm{d}x$ by a constant multiple of
\begin{align*}
\int_\R \bigg(\int_\R p_x(x,z)\,\mathrm{d}z\bigg)^2\,\frac{1}{m_1(x)}\,\mathrm{d}x
+ \int_\R \bigg(\int_\R f^2(z)\,p_x(x,z)\,\mathrm{d}z\bigg)^2\,\frac{1}{m_1(x)}\,\mathrm{d}x.
\end{align*}
By repeating the steps in \eqref{Dan's Jensen}, with $\psi=1$ and $\psi=f^2$, respectively, and with $(p,m_1)$ replacing $(\widetilde{p},\widetilde{m}_1)$, we see that the latter expression is finite. Hence, $\int_\R (G^{f^2;p})'(x)^2\,m_1(x)\,\mathrm{d}x<\infty$ and, due to the boundedness of $m_1$ away from $0$ on compact intervals (cf.\ \eqref{def:marginals}), the derivative $(G^{f^2;p})'$ is locally Lebesgue square integrable. Therefore, $G^{f^2;p}$ is locally $\frac{1}{2}$--H\"older continuous, as desired. \ep


\section{Proof of Theorem \ref{thm:uniq} and further discussion} \label{se:3}

In this section we first prove Theorem \ref{thm:uniq} and then give a more direct and enlightening argument for the independence of $X_t$ and $Y_t$ for each $t\ge0$ when $h \equiv f^2$.

\medskip

\noindent\textbf{Proof of Theorem \ref{thm:uniq}.} Recall from Theorem \ref{th:intro:main} that any stationary weak solution $(X_t,Y_t)_{t \ge 0}$ of the SDE \eqref{intro:mainSDEsystem} is strong and satisfies $X_t \stackrel{d}{=} \widehat X_t$ and $Y_t \stackrel{d}{=} \widehat Y_t$ for each $t \ge 0$, where $(\widehat X_t)_{t \ge 0}$ and $(\widehat Y_t)_{t \ge 0}$ are the unique stationary weak solutions of the SDEs \eqref{eq:Xhat} and \eqref{eq:Yhat}, respectively. Thanks to Proposition \ref{pr:general linear SDE} and strong uniqueness being a consequence of strong existence and weak uniqueness, it suffices to show that, for any stationary weak solution $(X_t,Y_t)_{t \ge 0}$ of \eqref{intro:mainSDEsystem}, $X_t$ and $Y_t$ are independent for each $t \ge 0$.

\medskip

Let $(X_t,Y_t)_{t \ge 0}$ be a stationary weak solution of \eqref{intro:mainSDEsystem}. By claim (iii) in Theorem \ref{th:intro:main}, the joint law ${\mathcal L}(X_0,Y_0)$ admits a density $p$. We introduce a transformation similar to that of Section \ref{se:2_2}. Define a new probability density $\widetilde p$ by
\begin{align*}
\widetilde p(x,y) = \frac{f^2(y)\,p(x,y)}{\int_{\R^2} f^2\,p\,\mathrm{d}x\,\mathrm{d}y}.
\end{align*}
We claim that $\widetilde p$ is the fixed-time marginal distribution of the stationary weak solution to
\begin{equation}\label{SDE_inde}
\begin{cases}
\;\mathrm{d}\widetilde{X}_t\,=\,b_1(\widetilde{X}_t)\,G^{f^2;p}(\widetilde{X}_t)\,\mathrm{d}t + \sigma_1(\widetilde{X}_t)\,\sqrt{G^{f^2;p}(\widetilde{X}_t)}\,\mathrm{d}W_t,
\\
\;\mathrm{d}\widetilde{Y}_t\;\,=\, b_2(\widetilde{Y}_t)\,f^{-2}(\widetilde{Y}_t)\,\mathrm{d}t+\sigma_2(\widetilde{Y}_t)\,f^{-1}(\widetilde{Y}_t)\,\mathrm{d}B_t.
\end{cases}
\end{equation}
To obtain this, note that $p$ itself solves \eqref{mainpde_weak}, which for $h \equiv f^2$ becomes
\begin{equation*}
\begin{split}
\forall \varphi \in C^\infty_c(\R^2):\quad
0 = \int_{\R^2}\bigg(\frac12\,\sigma_1^2(x)\, f^2(y)\,G^{f^2;p}(x)\,\varphi_{xx}(x,y) + \frac12\, \sigma_2^2(y)\,\varphi_{yy}(x,y) \qquad \qquad\; \\ 
+\,b_1(x)\,f^2(y)\,G^{f^2;p}(x)\,\varphi_x(x,y) + b_2(y)\,\varphi_y(x,y) 
\!\bigg) p(x,y)\,\mathrm{d}x\,\mathrm{d}y.
\end{split}
\end{equation*}
Equivalently,
\begin{equation*}
\begin{split}
\forall \varphi \in C^\infty_c(\R^2):\quad
0 = \int_{\R^2}\bigg(\frac12\,\sigma_1^2(x)\,G^{f^2;p}(x)\,\varphi_{xx}(x,y) + \frac12\,\sigma_2^2(y)\,f^{-2}(y)\,\varphi_{yy}(x,y) \qquad \qquad\; \\ 
+\,b_1(x)\,G^{f^2;p}(x)\,\varphi_x(x,y) + b_2(y)\,f^{-2}(y)\,\varphi_y(x,y) 
\!\bigg) \widetilde{p}(x,y)\,\mathrm{d}x\,\mathrm{d}y,
\end{split}
\end{equation*}
i.e., $\widetilde{p}$ solves the stationary Fokker-Planck equation associated with the SDE \eqref{SDE_inde}. Combining Proposition \ref{pr:general linear SDE} and \cite[Theorem 2.5]{trevisan} as in the beginning of the proof of Lemma \ref{lem:existence-estimate} we identify $\widetilde{p}$ as the fixed-time marginal distribution of the unique stationary weak solution to \eqref{SDE_inde}. The uniqueness of the latter renders $(\widetilde{X}_t)_{t\ge0}$ and $(\widetilde{Y}_t)_{t\ge0}$ independent, each being the unique stationary weak solution of the corresponding one-dimensional SDE (cf.\ Lemma \ref{le:intro:onedimensional}). Thus, $\widetilde p(x,y) = \widetilde{m}_1(x)\,\widetilde{m}_2(y)$, where $\widetilde{m}_1(x)=\int_\R \widetilde p(x,y)\,\mathrm{d}y$ and $\widetilde{m}_2(y)=\int_\R \widetilde p(x,y)\,\mathrm{d}x$. The definition of $\widetilde p$ implies that $p$ must also be of product form. \ep

\medskip

An interesting and more direct argument reveals why independent solutions arise when $h \equiv f^2$. We work formally here, implicitly assuming enough regularity for the differential equations to be valid in the classical sense, but the approach can be easily adapted to the setting of distributional solutions. Denote again by $m_1(x)$ and $m_2(y)$ the stationary probability densities for the one-dimensional marginal SDEs
\begin{align}
\mathrm{d}X_t = b_1(X_t)\,\mathrm{d}t + \sigma_1(X_t)\,\mathrm{d}W_t \quad\text{and}\quad \mathrm{d}Y_t = b_2(Y_t)\,\mathrm{d}t + \sigma_2(Y_t)\,\mathrm{d}B_t,	\label{intro:indepsolutions-1}
\end{align}
respectively. This means $m_1$ and $m_2$ solve the differential equations
\begin{align*}
\frac12 (\sigma_1^2 m_1)_{xx}-(b_1m_1)_x=0 \quad\text{and}\quad
\frac12 (\sigma_2^2 m_2)_{yy}-(b_2m_2)_y=0.
\end{align*}
Multiply the first equation by $m_2(y)f^2(y)$ and the second by $m_1(x)$, and then add the two resulting equations 
to find that $p(x,y)= m_1(x)\,m_2(y)$ solves the PDE
\begin{align*}
\frac12\,\big(\sigma_1^2(x)\,f^2(y)\,p(x,y)\big)_{xx}\!+\! \frac12\,\big(\sigma_2^2(x)\,p(x,y)\big)_{yy} 
\!-\!\big(b_1(x)\,f^2(y)\,p(x,y)\big)_x\!-\!\big(b_2(y)\,p(x,y)\big)_y\!=\!0.
\end{align*}
That is, $p(x,y)=m_1(x)\,m_2(y)$ is the stationary probability density for the SDE
\begin{align}
\begin{cases}
\;\mathrm{d}X_t = b_1(X_t)\,f^2(Y_t)\,\mathrm{d}t + \sigma_1(X_t)\,f(Y_t)\,\mathrm{d}W_t, \\
\;\mathrm{d}Y_t \;= b_2(Y_t)\,\mathrm{d}t + \sigma_2(Y_t)\,\mathrm{d}B_t.
\end{cases} \label{intro:indepsolutions-2}
\end{align}

\smallskip

Observe that it does not matter in the above argument whether we normalize $f^2(y)$ by the constant $\int_\R f^2\,m_2\,\mathrm{d}y$ or not, because the stationary probability density for the SDE
\[
\mathrm{d}X_t = a^2\,b_1(X_t)\,\mathrm{d}t + a\sigma_1(X_t)\,\mathrm{d}W_t
\]
is the same for each constant $a>0$. Indeed, the constant $a$ simply amounts to the time change $t \mapsto a^2t$ in the SDE. This point of view also lends to an intriguing interpretation of Theorem \ref{thm:uniq}. Given $(X_t)_{t\ge0}$ and $(Y_t)_{t \ge 0}$ solving the SDEs in \eqref{intro:indepsolutions-1} from independent initial conditions, suppose we define a time change by $t\mapsto\tau_t = \int_0^tf^2(Y_s)\,\mathrm{d}s$. Then $(X_{\tau_t},Y_t)_{t \ge 0}$ should be a weak solution of the SDE \eqref{intro:indepsolutions-2}. It is not immediately obvious, except when $f$ is constant, that the time-changed equation should admit the same (i.e., product form) stationary probability density, but Theorem \ref{thm:uniq} demonstrates that this is the case.

\bigskip\bigskip

\bibliographystyle{amsalpha}
\bibliography{Main}

\bigskip\bigskip\bigskip\bigskip

\end{document}